# Differentiation and Covering Constants for Hilbert-Schmidt and Quasi-Hilbert-Schmidt Operators

## Jinlu Li

Department of Mathematics
Shawnee State University
Portsmouth, OH 45662 USA
Email: jli@shawnee.edu

**Abstract** In this paper, we calculate the Fréchet derivatives and Mordukhovich derivatives (or coderivatives) of Hilbert-Schmidt operators on separable Hilbert spaces, by which we prove that the covering constant for Hilbert-Schmidt operators is zero. As an important class of Hilbert-Schmidt operators, we study the differentiability of Hilbert-Schmidt integral operators. Then, we introduce the concept of quasi-Hilbert-Schmidt operators on separable Hilbert spaces. We provide an example of quasi-Hilbert–Schmidt operators and find its Fréchet derivatives and Mordukhovich derivatives.

## 1. Introduction

The main theme of the theory of generalized differentiation in Banach spaces is deal with analysis of set-valued mappings, in which the cardinal topics primarily include the definitions of continuity and smoothness of considered set-valued mappings. It is well-known that Mordukhovich derivative (or coderivative) of set-valued mappings extends the concept of Fréchet differentiability for single-valued mappings in Banach spaces (See Theorem 1.38 in Mordukhovich [21]). Mordukhovich derivative provides a powerful concept for describing the generalized smoothness of set-valued mappings in Banach spaces, which has laid the foundation and has been playing the fundamental and crucial role in the theory of set-valued and variational analysis (See [21−24, 29]).

The theory of generalized differentiation (Mordukhovich derivative) in set-valued analysis has been potently developed and has been widely applied to many fields related set-valued mappings such as optimization, control theory, game theory, variational analysis, and so forth (see [1−7, 14−24, 29]). We believe that one of the most important applications of Mordukhovich derivatives is defining the covering constants for set-valued mappings in Banach spaces, which plays an important and crucial role in the Arutyunov Mordukhovich and Zhukovskiy Parameterized Coincidence Point Theorem (See [1, 21]). We list this theorem below for reader's easy reference.

**Theorem 3.1 in [1]** (Arutyunov Mordukhovich and Zhukovskiy Parameterized Coincidence Point Theorem). *Let $X$ and $Y$ be Asplund spaces and let $P$ be a topological space. Let $F(\cdot): X \rightrightarrows Y$ and $G(\cdot, \cdot): X \times P \rightrightarrows Y$ be set-valued mappings. Let $\bar{x} \in X$ and $\bar{y} \in Y$ with $\bar{y} \in F(\bar{x})$. Suppose that the following conditions are satisfied:*

(A1) *The multifunction $F: X \rightrightarrows Y$ is closed around $(\bar{x}, \bar{y})$.*

(A2) *There are neighborhoods $U \subset X$ of $\bar{x}$, $V \subset Y$ of $\bar{y}$, and $O$ of $\bar{p} \in P$ as well as a number $\beta \geq 0$ such*

*that the multifunction $G(\cdot, p): X \rightrightarrows Y$ is Lipschitz-like on U relative to V for each $p \in O$ with the uniform modulus $\beta$, while the multifunction $p \to G(\bar{x}, p)$ is lower/inner semicontinuous at $\bar{p}$.*

(A3) *The Lipschitzian modulus $\beta$ of $G(\cdot, p)$ is chosen as $\beta < \hat{\alpha}(F, \bar{x}, \bar{y})$, where $\hat{\alpha}(F, \bar{x}, \bar{y})$ is the covering constant of F around $(\bar{x}, \bar{y})$.*

*Then for each $\alpha > 0$ with $\beta < \alpha < \hat{\alpha}(F, \bar{x}, \bar{y})$, there exist a neighborhood $W \subset P$ of $\bar{p}$ and a single-valued mapping $\sigma: W \to X$ such that whenever $p \in W$ we have*

$$F(\sigma(p)) \cap G(\sigma(p), p) \neq \emptyset, \tag{1.1}$$

*and* $\qquad \|\sigma(p) - \bar{x}\|_X \leq \frac{\text{dist}(\bar{y}, G(\bar{x}, p))}{\alpha - \beta}. \tag{1.2}$

We see that the results of the Arutyunov Mordukhovich and Zhukovskiy Parameterized Coincidence Point Theorem are very strong and very useful. For example, it has been used for solving some parameterized optimal function problems in [1], stochastic fixed-point problems in [16], stochastic integral equations in [17], stochastic systems of linear equations in [18], stochastic vector variational inequalities and optimizations in [19] and parameterized equations in [20]. Based on the significant importance of the Arutyunov Mordukhovich and Zhukovskiy Parameterized Coincidence Point Theorem, it has been simply named as the AMZ Theorem (See [16−20]).

In order to use the AMZ Theorem to solve some real applications involved with set-valued mappings as appeared in Theorem 3.1 in [1], we have to precisely find the covering constant $\hat{\alpha}(F, \bar{x}, \bar{y})$, for the considered mapping $F$ at a given point $(\bar{x}, \bar{y})$. However, it is well-known that, this is the most difficult part in the real applications of the AMZ Theorem. Since the covering constant $\hat{\alpha}(F, \bar{x}, \bar{y})$ is defined by the Mordukhovich derivative $\hat{D}^*F(\bar{x}, \bar{y})$ of $F$ at point $(\bar{x}, \bar{y})$, more precisely speaking, the difficulty for calculating the covering constant $\hat{\alpha}(F, \bar{x}, \bar{y})$ partially comes from the computation of the Mordukhovich derivative $\hat{D}^*F(\bar{x}, \bar{y})$ of $F$ at point $(\bar{x}, \bar{y})$.

Hence, for a considered Mordukhovich differential mapping $F$ in Banach spaces, to find the exact Mordukhovich derivative $\hat{D}^*F(\bar{x}, \bar{y})$ of $F$ at some given point $(\bar{x}, \bar{y})$ will be a crucial topic in the applications of the AMZ Theorem, even though, in the theory of generalized differentiation.

For the above purposes, in this paper, we concentrate to study the covering constants for Hilbert-Schmidt operator on sparable Hilbert spaces, which includes Hilbert-Schmidt integral operator. We will extend the concepts of Hilbert-Schmidt operator to more general cases, which are called quasi-Hilbert-Schmidt operator. We will provide a concrete example of quasi-Hilbert-Schmidt operator on sparable Hilbert spaces and fined its Mordukhovich differentiability.

Notice that even for single-valued mappings, it is also very difficult to precisely find the Mordukhovich derivatives of the given mappings in Banach spaces, in general. Fortunately, the results of Theorem 1.38 in [21] provides close connections between the Fréchet derivative and the Mordukhovich derivative of single-valued mappings in Banach spaces. By using Theorem 1.38 in [21], the Mordukhovich derivative of Hilbert–Schmidt operators can be easily calculated, which comes from the Fréchet derivative of the considered Hilbert–Schmidt operators as linear and continuous single-valued mappings.

2. **Review for Fréchet Derivative and Mordukhovich Derivative (Coderivative) in Hilbert Spaces**

Mordukhovich derivatives (or coderivatives) of both single-valued and set-valued mappings are defined in Banach spaces. Since we concentrate to study covering constants for Hilbert-Schmidt Operators, which are considered in Hilbert spaces as underlying spaces, in this section, we briefly review the concepts and properties of Fréchet derivatives and Mordukhovich derivatives (or coderivatives) of single-valued mappings in Hilbert spaces. See [21–24, 29] for more details.

Let $(H, \|\cdot\|_H)$ be real Hilbert spaces with inner product $\langle \cdot, \cdot \rangle$ and origin $\theta$. For any $x \in H$ and $r > 0$, let $B(x, r)$ denote the closed ball in $H$ centered at point $x$ with radius $r$, respectively. Let $g: H \to H$ be a single-valued mapping and let $\bar{x} \in H$. If there is a linear and continuous mapping $\nabla g(\bar{x}): H \to H$ such that

$$\lim_{x \to \bar{x}} \frac{g(x) - g(\bar{x}) - \nabla g(\bar{x})(x - \bar{x})}{\|x - \bar{x}\|_H} = \theta,$$

then $g$ is said to be Fréchet differentiable at $\bar{x}$ and $\nabla g(\bar{x})$ is called the Fréchet derivative of $g$ at $\bar{x}$.

Then, we review the concept of Mordukhovich derivative (coderivative) for single-valued mappings in Hilbert spaces (See [21–24, 29] for more details). Let $\Delta$ be a nonempty subset of $H$ and let $g: \Delta \to H$ be a single-valued mapping. For $x \in \Delta$, a set-valued mapping $\widehat{D}^*g(x, g(x)): H \to H$ is defined, for any $y \in H$, by (See Definitions 1.13 and 1.32 in Chapter 1 in [21])

$$\widehat{D}^*g(x, g(x))(y) = \left\{ z \in H: \limsup_{\substack{(u, g(u)) \to (x, g(x)) \\ u \in \Delta}} \frac{\langle z, u - x \rangle - \langle y, g(u) - g(x) \rangle}{\|u - x\|_H + \|g(u) - g(x)\|_H} \leq 0 \right\}.$$

If this set-valued mapping $\widehat{D}^*g(x, g(x)): H \to H$ satisfies that

$$\widehat{D}^*g(x, g(x))(y) \neq \emptyset, \text{ for every } y \in H, \tag{2.1}$$

then $g$ is said to be Mordukhovich differentiable (codifferentiable) at the point $x$ and $\widehat{D}^*g(x, g(x))$ is called the Mordukhovich derivative (which is also called Mordukhovich coderivative, or coderivative) of $g$ at point $x$. For this single-valued mapping $g: \Delta \to H$, at point $x$, we write

$$\widehat{D}^*g(x, g(x))(y) \equiv \widehat{D}^*g(x)(y), \text{ for every } y \in H.$$

Furthermore, if the single-valued mapping $g: \Delta \to H$ is a continuous mapping, then, by the above definition, the Mordukhovich derivative of $g$ at point $x$ is calculated by

$$\widehat{D}^*g(x)(y) = \left\{ z \in H: \limsup_{\substack{u \to x \\ u \in \Delta}} \frac{\langle z, u - x \rangle - \langle y, g(u) - g(x) \rangle}{\|u - x\|_H + \|g(u) - g(x)\|_H} \leq 0 \right\}, \text{ for any } y \in H. \tag{2.2}$$

The following theorem shows the connection between Fréchet derivatives and Mordukhovich derivatives for sing-valued mappings. The results of the following theorem provide a powerful tool to calculate the Mordukhovich derivatives by the Fréchet derivatives of single-valued mappings.

**Theorem 1.38 in [21] in Hilbert spaces**. *Let $H$ be a Hilbert and let $g: H \to H$ be a single-valued mapping. Suppose that $g$ is Fréchet differentiable at $x \in H$. Then, the Mordukhovich derivative of $g$ at $x$ satisfies the following equation*

$$\widehat{D}^*g(x)(y) = \{(\nabla g(x))^*(y)\}, \text{ for all } y \in H.$$

One of the important applications of Mordukhovich derivatives of set-valued mappings is to define the covering constants for set-valued mappings. The covering constant for $\Phi: H \rightrightarrows H$ at point $(\bar{x}, \bar{y}) \in \text{gph } \Phi$ is defined by (see (2.4) in [1], also see [21])

$$\hat{\alpha}(\Phi,\bar{x},\bar{y}) := \sup_{\eta>0} \inf\{\|z\|_H : z \in \widehat{D}^*\Phi(x,y)(w), x \in \mathbb{B}(\bar{x},\eta), y \in \Phi(x) \cap \mathbb{B}(\bar{y},\eta), \|w\|_H = 1\}. \quad (2.4)$$

Let $\mathbb{B}$ denote the unit closed ball in $H$. For $\bar{x}, \bar{y} \in H$ and $\eta > 0$, let $\mathbb{B}(\bar{x},\eta)$ denote the closed ball in $H$ centered at $\bar{x}$ with radius $\eta$, and $\mathbb{B}(\bar{y},\eta)$ denote the closed ball in $H$ centered at $\bar{y}$ with radius $\eta$. In particular, let $g: X \to Y$ be a single-valued mapping. For any $\bar{x}, \bar{y} \in H$ with $\bar{y} = g(\bar{x})$, (2.4) becomes

$$\hat{\alpha}(g,\bar{x},\bar{y}) = \sup_{\eta>0} \inf\{\|z\|_H : z \in \widehat{D}^*g(x)(w), x \in \mathbb{B}(\bar{x},\eta), g(x) \in \mathbb{B}(\bar{y},\eta), \|w\|_H = 1\}. \quad (2.5)$$

**Definition 1.51 in [21] in Hilbert spaces** (covering properties) Let $F: H \rightrightarrows H$ with dom $F \neq \emptyset$.

(i) Let $U$ and $V$ be nonempty subsets in $H$. We say that $F$ enjoys the covering property on $U$ relative to $V$ if there is $\gamma > 0$ such that

$$F(x) \cap V + \gamma r \mathbb{B} \subset F(x + r\mathbb{B}), \text{ whenever } x + r\mathbb{B} \subset U, \text{ as } r > 0. \quad (2.2)$$

(ii) Given $(\bar{x},\bar{y}) \in \operatorname{gph} F$, we say that $F$ has the local covering property around $(\bar{x},\bar{y})$ with modulus $\gamma > 0$ if there is a neighborhood $U$ of $\bar{x}$ and a neighborhood $V$ of $\bar{y}$ such that the above inclusion in (i) holds. The supremum of all such moduli $\{\gamma\}$ is called the exact covering bound of $F$ around $(\bar{x},\bar{y})$, which is denoted by

$$\operatorname{cov} F(\bar{x},\bar{y}) = \sup\{\gamma : \gamma \text{ satisfies (2.2) for some } U, V \subset H\}.$$

(iii) (Locally covering properties for single-valued mappings (2.1) in [1]) Let $x_0 \in X$ and $g: H \to H$ be a single-valued mapping that is assumed to be continuous in a neighborhood of $x_0$ (or, more general, assumed to be defined in a neighborhood of $x_0$). $g$ is said to be locally $\alpha$-covering in a neighborhood of $x_0$ if there is a $\delta > 0$ such that

$$\mathbb{B}(x,r) \subset \mathbb{B}(x_0,\delta) \implies \mathbb{B}(g(x), \alpha r) \subset g(\mathbb{B}(x,r)), \text{ for any } r > 0.$$

If there is $\alpha > 0$ such that $g$ is locally $\alpha$-covering in a neighborhood of $x_0$, then, we write

$$\alpha(g, x_0) = \sup\{\alpha > 0 : F \text{ is locally } \alpha - \text{covering in a neighborhood of } x_0\}.$$

Otherwise, we denote $\alpha(g, x_0) = 0$.

Since every Hilbert space is a Asplund space (See [9, 11-13, 27]), then, Theorem 4.1, Corollary 4.3 in [21] and Theorem 2.1 in [1] all hold for Hilbert spaces.

**Theorem 4.1 in [21]** (neighborhood characterization of local covering) *Let $F: H \rightrightarrows H$ be a set-valued mapping. Assume that $F$ is closed-graph around $(\bar{x},\bar{y}) \in \operatorname{gph} F$. Then the following are equivalent*:

(a) *$F$ enjoys the local covering property around $(\bar{x},\bar{y})$ (that is, $\operatorname{cov} F(\bar{x},\bar{y}) > 0$).*
(b) *One has $\hat{\alpha}(F,\bar{x},\bar{y}) > 0$.*

*Moreover, the exact covering bound of $F$ around $(\bar{x},\bar{y})$ is computed by*

$$\operatorname{cov} F(\bar{x},\bar{y}) = \hat{\alpha}(F,\bar{x},\bar{y}).$$

**Corollary 4.3 in [21] in Hilbert spaces** (neighborhood covering criterion for single-valued mappings) *Let $g: H \to H$ be a single-valued mapping. Assume that $g$ is Lipschitz continuous around some point $\bar{x}$. Then the conclusions of Theorem 4.1 hold with the covering constant $\hat{\alpha}(g,\bar{x})$ computed by*

$$\mathrm{cov} g(\bar{x}) = \hat{\alpha}(g, \bar{x})$$

**Theorem 2.1 in [1] in Hilbert spaces** (The covering criterion). *Let $\Phi: H \rightrightarrows H$ be a set-valued mapping. Let $(\bar{x}, \bar{y}) \in \mathrm{gph}\,\Phi$, and let $\hat{\alpha}(\Phi, \bar{x}, \bar{y})$ be taken from (2.5). Then:*

(i) *The α-covering property of $\Phi$ around $(\bar{x}, \bar{y})$ implies that $\alpha \leq \hat{\alpha}(\Phi, \bar{x}, \bar{y})$.*

(ii) *If $0 < \alpha < \hat{\alpha}(\Phi, \bar{x}, \bar{y})$, then $\Phi$ has the α-covering property around $(\bar{x}, \bar{y})$ provided that its graph is closed around this point.*

## 3. Hilbert–Schmidt Operator

In this section, we review the concepts and properties of Hilbert-Schmidt operator on Hilbert spaces. See [8, 10, 25–26, 28, 30–32] for more details.

In this paper, unless otherwise stated, let $(H, \|\cdot\|_H)$ be a separable Hilbert space with inner product $\langle \cdot, \cdot \rangle$ and equipped with an orthonormal (Schauder) basis $\{e_i: i \in \mathbb{N}\}$. Based on this basis, the elements of $H$ have the following representation.

$$x = \sum_{i=1}^{\infty} \langle x, e_i \rangle e_i, \text{ for any } x \in H.$$

The sequence of coefficients of $x$ with respect to the orthonormal basis $\{e_i: i \in \mathbb{N}\}$ is denoted by

$$\mathrm{sc}(x) = \{\langle x, e_i \rangle : i \in \mathbb{N}\} = (\langle x, e_1 \rangle, \langle x, e_2 \rangle, \dots) \in l_2, \text{ for } x \in H.$$

Let $T: H \to H$ be a linear and continuous (bonded) operator. If $T$ satisfies that

$$\sum_{i=1}^{\infty} \|Te_i\|_H^2 < \infty. \tag{3.1}$$

then, $T$ is called a Hilbert-Schmidt operator on $H$. The sum in (3.1) defines the Hilbert-Schmidt norm $\|T\|_{HS}$ (or the Frobenius norm if $H$ is Euclidean space) of $T$ such that

$$\|T\|_{HS}^2 = \sum_{i=1}^{\infty} \|Te_i\|_H^2 = \sum_{i,j=1}^{\infty} \langle Te_i, e_j \rangle^2. \tag{3.2}$$

We list below some properties of Hilbert-Schmidt operators, which will be used in this paper (See Conway [10]).

(HS1) Every Hilbert–Schmidt operator $T: H \to H$ is a compact operator.

(HS2) If $T: H \to H$ is a Hilbert–Schmidt operator, then we have $\|T\| \leq \|T\|_{HS}$.

In this section, we use the property (HS1) for Hilbert-Schmidt operators to prove that, for any Hilbert-Schmidt operator $T$ in a Hilbert space $(H, \|\cdot\|_H)$, if the considered Hilbert space $H$ has infinity dimension, then, for any given nonempty subsets $U$ and $V$ in $H$, $T$ does not enjoy the covering property on $U$ relative to $V$. This implies that, for any $\bar{x}, \bar{y} \in H$, with $\bar{y} = T\bar{x}$, the covering constant for $T$ at point $(\bar{x}, \bar{y})$ is zero.

**Theorem 3.1**. *Suppose that $H$ is an infinity dimensional Hilbert space. Then, for every Hilbert-Schmidt operator $T: H \to H$, for any $\bar{x}, \bar{y} \in H$, with $\bar{y} = T\bar{x}$, we have that*

(i) *For any given nonempty subsets $U$ and $V$ in $H$, $T$ does not enjoy the covering property on $U$ relative to $V$, which implies that the exact covering bound of $T$ around $(\bar{x}, \bar{y})$ is zero; that is*

$$\mathrm{cov} T(\bar{x}, \bar{y}) = 0.$$

(ii)    *The covering constant of T is constant with*

$$\hat{\alpha}(T, \bar{x}, \bar{y}) = 0.$$

*Proof.* By the property (HS1) of Hilbert-Schmidt operators on $H$, operator $T: H \to H$ is a compact operator. For any given nonempty subsets $U$ and $V$ in $H$, and for any $r > 0$, assume that $\bar{x} + r\mathbb{B} \subset U$. Since $\bar{x} + r\mathbb{B}_H$ is a closed and bounded subset of $H$, then $T(\bar{x} + r\mathbb{B})$ is a restively compact subset of $H$. On the other hand, for any $\gamma > 0$, by $\bar{y} = T\bar{x}$, we have

$$T(\bar{x}) \cap V + \gamma r\mathbb{B} = \begin{cases} \bar{y} + \gamma r\mathbb{B}, & \text{if } \bar{y} \in V, \\ \gamma r\mathbb{B}, & \text{if } \bar{y} \notin V. \end{cases}$$

This implies that the set $T(\bar{x}) \cap V + \gamma r\mathbb{B}$ is a closed ball with radius $\gamma r$ centered at either $\bar{y}$ or $\theta$. Since $H$ is an infinity dimensional Hilbert space, this closed ball $T(\bar{x}) \cap V + \gamma r\mathbb{B}$ cannot be contained in a restively compact subset of $H$. That is,

$$T(\bar{x}) \cap V + \gamma r\mathbb{B} \nsubseteq T(\bar{x} + r\mathbb{B}), \text{ if } \bar{x} + r\mathbb{B} \subset U.$$

This implies that

$$\operatorname{cov} T(\bar{x}, \bar{y}) = 0, \text{ for any } \bar{x}, \bar{y} \in H, \text{ with } \bar{y} = T\bar{x}. \tag{3.3}$$

Since $T: H \to H$ is a linear and continuous (Hilbert-Schmidt) operator and $T$ is a compact operator, this implies that $T$ is closed-graph around an arbitrarily given point $(\bar{x}, \bar{y}) \in \operatorname{gph} T$. Notice that every Hilbert space is a Asplund spaces. Then, by Theorem 4.1 in [21] (neighborhood characterization of local covering on Asplund spaces) and by (3.3), we have

$$\hat{\alpha}(F, \bar{x}, \bar{y}) = \operatorname{cov} T(\bar{x}, \bar{y}) = 0, \text{ for any } \bar{x}, \bar{y} \in H, \text{ with } \bar{y} = T\bar{x}. \qquad \square$$

Recall that in functional analysis, as a branch of mathematics, a compact operator is a linear operator. However, we see that, in the proof of Theorem 3.1, the compactness of the considered Hilbert-Schmidt operator $T$ plays a crucial important roll. Meanwhile, the linearity of $T$ is not used in the proof; and therefore, the linearity of $T$ is not required for Theorem 3.1. Hence, with the similar proof, Theorem 3.1 can be extended to the following more general theorem.

**Theorem 3.2**. *Suppose that $B$ is an infinity dimensional Banach space. Let $S: B \to B$ be a mapping. Let $\bar{x}, \bar{y} \in B$, with $\bar{y} = S\bar{x}$. Suppose that $S$ satisfies the following conditions at point $(\bar{x}, \bar{y})$,*

(C1) *$S$ is closed-graph around point $(\bar{x}, \bar{y}) \in \operatorname{gph} S$;*
(C2) *$S$ maps every bounded subset of $B$ to a relatively compact subset of $B$.*

*Then, we have*

(i)    *For any given nonempty subsets $U$ and $V$ in $B$, $S$ does not enjoy the covering property on $U$ relative to $V$, which implies that the exact covering bound of $S$ around $(\bar{x}, \bar{y})$ is zero; that is*

$$\operatorname{cov} S(\bar{x}, \bar{y}) = 0.$$

(ii)   *The covering constant of $S$ is constant with*

$$\hat{\alpha}(S, \bar{x}, \bar{y}) = 0.$$

*Proof.* The proof of this theorem is similar to the proof of Theorem 3.1, which is omitted here. □

### 4. An Example of Hilbert–Schmidt Operator in Separable Hilbert Spaces

Recall that, as we stated in the previous section, $(H, \|\cdot\|_H)$ is a Hilbert space with inner product $\langle \cdot, \cdot \rangle$ and with an orthonormal (Schauder) basis $\{e_i : i \in \mathbb{N}\}$. In this section, we consider a special Hilbert-Schmidt operator on $H$, which is defined by a double sequence. We investigate its Fréchet and Mordukhovich derivatives. By Theorem 3.1, its covering constant is zero.

Let $A = \{a_{ij} : i, j \in \mathbb{N}\}$ be a double sequence of real numbers. Suppose that $A$ satisfies that following conditions:

$$\sum_{i,j} a_{ij}^2 < \infty. \tag{4.1}$$

Associating to this double sequence $A$, we define a mapping $T: H \to H$, for any $x \in H$, by

$$Tx = (\langle x, e_1 \rangle, \langle x, e_2 \rangle, \ldots) A \begin{pmatrix} e_1 \\ e_2 \\ \vdots \end{pmatrix} = \mathrm{sc}(x) A \begin{pmatrix} e_1 \\ e_2 \\ \vdots \end{pmatrix}. \tag{4.2}$$

**Theorem 4.1.** *Let $T: H \to H$ be the operator defined by* (4.2) *associated to the double matrix $A$ satisfying condition* (4.1). *Then $T$ has the following properties.*

(i) *$T: H \to H$ is a linear and continuous Hilbert-Schmidt operator on $H$. The Hilbert–Schmidt norm of $T$ satisfies that*

$$\|T\|_{HS}^2 = \sum_{i=1}^{\infty} \sum_{j=1}^{\infty} a_{ij}^2.$$

(ii) *$T: H \to H$ is Fréchet differentiable on $H$ such that*

$$\nabla g(z) = T, \textit{for any } z \in H;$$

(iii) *$T: H \to H$ is Mordukhovich differentiable on $H$ such that*

$$\widehat{D}^* g(z) = T^*, \textit{for any } z \in H.$$

*Here, $T^*$ is the adjoint operator of $T$, which is defined by $A^T$.*

(iv) *The covering constant of $T$ is constant with*

$$\hat{\alpha}(T, \bar{x}, \bar{y}) = 0, \bar{x}, \bar{y} \in H, \textit{with } \bar{y} = H\bar{x}.$$

*Proof.* By definition (3.4), for any $x \in H$, we have

$$Tx = \mathrm{sc}(x) A \begin{pmatrix} e_1 \\ e_2 \\ \vdots \end{pmatrix} = (\langle x, e_1 \rangle, \langle x, e_2 \rangle, \ldots) A \begin{pmatrix} e_1 \\ e_2 \\ \vdots \end{pmatrix}$$

$$= \left( \sum_{i=1}^{\infty} a_{i1} \langle x, e_i \rangle, \sum_{i=1}^{\infty} a_{i2} \langle x, e_i \rangle, \ldots \right) \begin{pmatrix} e_1 \\ e_2 \\ \vdots \end{pmatrix}$$

$$= \sum_{j=1}^{\infty} \sum_{i=1}^{\infty} a_{ij} \langle x, e_i \rangle e_j. \tag{4.3}$$

In particular, for any $k \in \mathbb{N}$, we have

$$Te_k == \text{sc}(e_k) A \begin{pmatrix} e_1 \\ e_2 \\ \vdots \end{pmatrix} = (\langle e_k, e_1 \rangle, \langle e_k, e_2 \rangle, \ldots) A \begin{pmatrix} e_1 \\ e_2 \\ \vdots \end{pmatrix}$$

$$= \sum_{j=1}^{\infty} \sum_{i=1}^{\infty} a_{ij} \langle e_k, e_i \rangle e_j$$

$$= \sum_{j=1}^{\infty} a_{kj} e_j. \tag{4.4}$$

It is clear that $T: H \to H$ is linear. For any $x \in H$, we calculate

$$\|Tx\|_H^2 = \sum_{j=1}^{\infty} \left( \sum_{i=1}^{\infty} a_{ij} \langle x, e_i \rangle \right)^2$$

$$\leq \sum_{j=1}^{\infty} \sum_{i=1}^{\infty} a_{ij}^2 \sum_{i=1}^{\infty} \langle x, e_i \rangle^2$$

$$= \sum_{i,j} a_{ij}^2 \|x\|_H^2. \tag{4.5}$$

By the assumption (4.1), this implies that $T: H \to H$ is a linear and continuous mapping. (4.5) implies that the norm of the mapping $T$ satisfies that

$$\|T\| \leq \sqrt{\sum_{i,j} a_{ij}^2}.$$

Similarly to (4.4), by (4.3) and (4.1), we check that

$$\|T\|_{HS}^2 = \sum_{k=1}^{\infty} \|Te_k\|_H^2 = \sum_{k=1}^{\infty} \left\| \sum_{j=1}^{\infty} a_{kj} e_j \right\|_H^2 = \sum_{k=1}^{\infty} \sum_{j=1}^{\infty} a_{kj}^2 < \infty.$$

This implies that $T: H \to H$ is a linear and continuous Hilbert-Schmidt operator on $H$. The Hilbert–Schmidt norm of $T$ satisfies that

$$\|T\|_{HS}^2 = \sum_{i=1}^{\infty} \sum_{j=1}^{\infty} a_{ij}^2. \tag{3.7}$$

This proves (i). Then, we prove (ii) and (iii). Since that $T: H \to H$ is a linear and continuous on $H$, we obtain that $T$ is Fréchet differentiable on $H$ such that

$$\nabla g(z) = T, \text{ for any } z \in H.$$

This proves (ii). By Theorem 1.38 in [21], $T$ is Mordukhovich differentiable on $H$ such that

$$\widehat{D}^* g(z) = T^*, \text{ for any } z \in H.$$

Here, $T^*$ is the adjoint operator of $T$, which is defined by $A^T$. This proves (iii). We finally prove (iv). Notice that part (iv) follows from Theorem 3.1 immediately. However, for the purposes of deeply understanding the concept of covering constants, we spend some time to directly prove the results of (iv) by using the Definition (2.7). Similar to (4.3), for any $z, y \in H$, we have

$$\widehat{D}^* T(z)(y) = T^* y = \text{sc}(y) A^T \begin{pmatrix} e_1 \\ e_2 \\ \vdots \end{pmatrix} = (\langle y, e_1 \rangle, \langle y, e_2 \rangle, \ldots) A^T \begin{pmatrix} e_1 \\ e_2 \\ \vdots \end{pmatrix}$$

$$= \left(\sum_{j=1}^{\infty} a_{1j} \langle y, e_j \rangle, \sum_{j=1}^{\infty} a_{2j} \langle y, e_j \rangle, \ldots \right) \begin{pmatrix} e_1 \\ e_2 \\ \vdots \end{pmatrix}$$

$$= \sum_{i=1}^{\infty} \left(\sum_{j=1}^{\infty} a_{ij} \langle y, e_j \rangle \right) e_i. \tag{4.6}$$

In particular, for any $k = 1, 2, \ldots$, we have

$$T^* e_k = \mathrm{sc}(e_k) A^T \begin{pmatrix} e_1 \\ e_2 \\ \vdots \end{pmatrix} = (\langle e_k, e_1 \rangle, \langle e_k, e_2 \rangle, \ldots ) A^T \begin{pmatrix} e_1 \\ e_2 \\ \vdots \end{pmatrix}$$

$$= \sum_{i=1}^{\infty} \left(\sum_{j=1}^{\infty} a_{ij} \langle e_k, e_j \rangle \right) e_i$$

$$= \sum_{i=1}^{\infty} a_{ik} e_i. \tag{4.7}$$

For any $x, y \in H$, by (4.6) and (4.7), we calculate

$$\langle Tx, y \rangle = \langle \sum_{j=1}^{\infty} \sum_{i=1}^{\infty} a_{ij} \langle x, e_i \rangle e_j, \sum_{j=1}^{\infty} \langle y, e_j \rangle e_j \rangle$$

$$= \sum_{j=1}^{\infty} \sum_{i=1}^{\infty} a_{ij} \langle x, e_i \rangle \langle y, e_j \rangle$$

$$= \sum_{i=1}^{\infty} \sum_{j=1}^{\infty} a_{ij} \langle y, e_j \rangle \langle x, e_i \rangle$$

$$= \langle \sum_{i=1}^{\infty} \langle x, e_i \rangle e_i, \sum_{i=1}^{\infty} \left(\sum_{j=1}^{\infty} a_{ij} \langle y, e_j \rangle \right) e_i \rangle$$

$$= \langle x, T^* y \rangle.$$

By the assumption (4.1), we have

$$\sum_{k=1}^{\infty} \sum_{i=1}^{\infty} a_{ik}^2 = \sum_{i,j} a_{ij}^2 < \infty.$$

This implies that

$$\sum_{i=1}^{\infty} a_{ik}^2 \to 0, \text{ as } k \to \infty. \tag{4.8}$$

By (2.7), for any $\bar{x}, \bar{y} \in H$ with $\bar{y} = T\bar{x}$, by (4.6), (4.7) and (4.8), we calculate the covering constant.

$$\hat{\alpha}(T, \bar{x}, \bar{y}) = \sup_{\eta > 0} \inf \{ \|z\|_H : z \in \widehat{D}^* T(x)(y), x \in \mathbb{B}(\bar{x}, \eta), T(x) \in \mathbb{B}(\bar{y}, \eta), \|y\|_H = 1 \}$$

$$= \sup_{\eta > 0} \inf \{ \|T^* y\|_H : T^* y, x \in \mathbb{B}(\bar{x}, \eta), T(x) \in \mathbb{B}(\bar{y}, \eta), \|y\|_H = 1 \}$$

$$\leq \sup_{\eta > 0} \inf \{ \|\sum_{i=1}^{\infty} a_{ik} e_i\|_H : T^* y, x \in \mathbb{B}(\bar{x}, \eta), T(x) \in \mathbb{B}(\bar{y}, \eta), \|y\|_H = 1, y = e_k \}$$

$$= \sup_{\eta > 0} \inf \left\{ \sqrt{\sum_{i=1}^{\infty} a_{ik}^2} : T^* y, x \in \mathbb{B}(\bar{x}, \eta), T(x) \in \mathbb{B}(\bar{y}, \eta), \|y\|_H = 1, y = e_k \right\}$$

$$= \lim_{k \to \infty} \sqrt{\sum_{i=1}^{\infty} a_{ik}^2}$$

$$= 0.$$

This is that

$$\hat{\alpha}(T, \bar{x}, \bar{y}) = 0, \text{ for any } \bar{x}, \bar{y} \in H \text{ with } \bar{y} = T\bar{x}. \qquad \square$$

## 5. Hilbert-Schmidt Integral Operator

As an important class of Hilbert–Schmidt operators, in this section, we consider Hilbert-Schmidt integral operators. We first review the concepts of Hilbert–Schmidt integral operators. For more details, see Bump [8], Renardy and Rogers [28] and Simon [31]. Let $n \geq 1$ and let $(\mathbb{R}^n, \|\cdot\|_n)$ be the standard $n$-d Euclidean space with the ordinal Hilbert $L_2$-norm denoted by $\|\cdot\|_n$ and with row vectors. Let $\theta$ denote the origin of $\mathbb{R}^n$. Let $\Omega$ be a domain in $\mathbb{R}^n$ ($\Omega$ is a nonempty, connected and open subset of $\mathbb{R}^n$). Let $(L^2(\Omega), \|\cdot\|_{L^2(\Omega)})$ and $(L^2(\Omega \times \Omega), \|\cdot\|_{L^2(\Omega \times \Omega)})$ denote the corresponding standard Hilbert spaces of real valued functions on $\Omega$ and $\Omega \times \Omega$, respectively. Let $\langle \cdot, \cdot \rangle$ denote the inner product on $L^2(\Omega)$. Let $k: \Omega \times \Omega \to \mathbb{R}$ be a real valued function such that

$$\|k\|^2_{L^2(\Omega \times \Omega)} = \int_\Omega \int_\Omega |k(s,t)|^2 ds dt < \infty. \tag{5.1}$$

Then, $k$ is called a Hilbert-Schmidt kernel. The associated integral operator $T: L^2(\Omega) \to L^2(\Omega)$ is defined, for any $f \in L^2(\Omega)$, by

$$T(f)(s) = \int_\Omega k(s,t) f(t) dt, \text{ for any } s \in \Omega. \tag{5.2}$$

It is well-known that $T: L^2(\Omega) \to L^2(\Omega)$ is a (linear and continuous) Hilbert-Schmidt integral operator on $L^2(\Omega)$ with the Hilbert-Schmidt kernel $k$. The Hilbert–Schmidt norm of this Hilbert-Schmidt integral operator $T: L^2(\Omega) \to L^2(\Omega)$ satisfies that (See Conway [10], Renardy and Rogers [28] and Simon [31]).

$$\|T\|_{HS} = \|k\|_{L^2(\Omega \times \Omega)} = \sqrt{\int_\Omega \int_\Omega |k(s,t)|^2 ds dt}. \tag{5.3}$$

Let $\{\varphi_i : i \in \mathbb{N}\}$ be an orthonormal (Schauder) basis of $L^2(\Omega)$. By Fubini's theorem and the $\sigma$-finiteness of the orthonormal basis for $L^2(\Omega)$, $\{\varphi_i \varphi_j : i, j \in \mathbb{N}\}$ is an orthonormal basis for $L^2(\Omega \times \Omega)$. Thus, for the considered kernel $k$, the double sequence of coefficients of $k$ is written as $\{c_{ij} : i, j \in \mathbb{N}\}$, which is defined, for each $i, j \in \mathbb{N}$, by

$$c_{ij} = \int_{\Omega \times \Omega} k(s,t) \varphi_i(s) \varphi_j(t) ds dt = \int_\Omega \varphi_j(t) \left( \int_\Omega k(s,t) \varphi_i(s) ds \right) dt,$$

then, the Hilbert-Schmidt kernel $k$ of $T$ has the following representation, in the sense of $L^2(\Omega \times \Omega)$,

$$k(s,t) = \sum_{i,j=1}^\infty c_{ij} \varphi_i(s) \varphi_j(t), \text{ for every } (s,t) \in \Omega \times \Omega. \tag{5.4}$$

By the square integrability of $k$ on $\Omega \times \Omega$, with respect to the orthonormal basis $\{\varphi_i \varphi_j : i, j \in \mathbb{N}\}$ for $L^2(\Omega \times \Omega)$, and by (5.4), we have

$$\|k\|^2_{L^2(\Omega \times \Omega)} = \int_\Omega \int_\Omega |k(s,t)|^2 ds dt$$

$$= \int_\Omega \int_\Omega \left| \sum_{i,j=1}^\infty c_{ij} \varphi_j(t) \varphi_i(s) \right|^2 ds dt$$

$$= \sum_{i,j=1}^{\infty} c_{ij}^2.$$

By (5.3) and (5.1), this implies that the Hilbert-Schmidt norm of $T$ can be calculated by

$$\|T\|_{HS}^2 = \|k\|_{L^2(\Omega \times \Omega)}^2 = \sum_{i,j=1}^{\infty} c_{ij}^2 < \infty.$$

Meanwhile, the operator $T: L^2(\Omega) \to L^2(\Omega)$ defined in (5.2) can be rewritten, for any $f \in L^2(\Omega)$, by

$$T(f)(s) = \sum_{i,j=1}^{\infty} c_{ij} \left( \int_{\Omega} \varphi_j(t) f(t) dt \right) \varphi_i(s)$$

$$= \sum_{i,j=1}^{\infty} c_{ij} \langle f, \varphi_j \rangle \varphi_i(s)$$

$$= \sum_{i=1}^{\infty} \left( \sum_{j=1}^{\infty} c_{ij} \langle f, \varphi_j \rangle \right) \varphi_i(s), \text{ for any } s \in \Omega. \tag{5.5}$$

**Theorem 5.1**. *The Hilbert–Schmidt integral operator $T: L^2(\Omega) \to L^2(\Omega)$ defined by (5.2) has the following differentiational properties.*

(i)      *$T$ is Fréchet differentiable on $L^2(\Omega)$ such that, for any $h \in L^2(\Omega)$, we have*

$$\nabla T(h)(f) = Tf = \sum_{i=1}^{\infty} \left( \sum_{j=1}^{\infty} c_{ij} \langle f, \varphi_j \rangle \right) \varphi_i, \text{ for any } f \in L^2(\Omega).$$

(ii)     *$T$ is Mordukhovich differentiable on $L^2(\Omega)$ such that, for any $h \in L^2(\Omega)$, we have*

$$\widehat{D}^* T(h)(g) = T^* g = \sum_{j=1}^{\infty} \left( \sum_{i=1}^{\infty} c_{ij} \langle g, \varphi_i \rangle \right) \varphi_j(t), \text{ for any } g \in L^2(\Omega).$$

(iii)    *The covering constant of $T$ is constant with*

$$\hat{\alpha}(T, f, u) = 0, \text{ for any } f, u \in L^2(\Omega) \text{ with } u = Tf.$$

*Proof.* Proof of (i). Since $T: L^2(\Omega) \to L^2(\Omega)$ is a linear and continuous mapping, then $T$ is Fréchet differentiable on $L^2(\Omega)$ such that, for any $h \in L^2(\Omega)$, by (5.5), the Fréchet derivative of $T$ at $h$ satisfies that, for any $f \in L^2(\Omega)$,

$$\nabla T(h)(f)(s) = Tf(s)$$

$$= \sum_{i=1}^{\infty} \left( \sum_{j=1}^{\infty} c_{ij} \langle f, \varphi_j \rangle \right) \varphi_i(s), \text{ for any } s \in \Omega.$$

Proof of (ii). The adjoint operator of $T$ is denoted by $T^*: L^2(\Omega) \to L^2(\Omega)$ such that, for any $g \in L^2(\Omega)$,

$$T^*(g)(t) = \int_{\Omega} k(s,t) g(s) ds$$

$$= \sum_{i,j=1}^{\infty} c_{ij} \left( \int_{\Omega} \varphi_i(s) g(s) dt \right) \varphi_j(t)$$

$$= \sum_{i,j=1}^{\infty} c_{ij} \langle g, \varphi_i \rangle \varphi_j(t)$$

$$= \sum_{j=1}^{\infty} \left( \sum_{i=1}^{\infty} c_{ij} \langle g, \varphi_i \rangle \right) \varphi_j(t), \text{ for any } t \in \Omega. \tag{5.6}$$

By (5.2), (5.3) and (5.4), and by Fubini's theorem, one has

$$\langle T(f), g \rangle = \langle \int_\Omega k(s,t)f(t)dt, g \rangle$$

$$= \int_\Omega \left( \int_\Omega k(s,t)f(t)dt \right) g(s)ds$$

$$= \int_\Omega \int_\Omega k(s,t)\, g(s)f(t)dsdt$$

$$= \int_\Omega \left( \int_\Omega k(s,t)\, g(s)ds \right) f(t)dt$$

$$= \langle f, T^*(g) \rangle.$$

This proves that $T^*: L^2(\Omega) \to L^2(\Omega)$ defined by (5.6) is the adjoint operator of $T$ satisfying

$$\langle T(f), g \rangle = \langle f, T^*(g) \rangle, \text{ for any } f, g \in L^2(\Omega).$$

By Definition (5.6), $T^*: L^2(\Omega) \to L^2(\Omega)$ has the same properties with $T: L^2(\Omega) \to L^2(\Omega)$.

(a) $T^*$ is also a Hilbert-Schmidt integral operator with the same Hilbert–Schmidt norm

$$\|T^*\|_{HS} = \|T\|_{HS} = \|k\|_{L^2(\Omega \times \Omega)} < \infty.$$

(b) $T^*$ is also both continuous and compact linear operator.

By Theorem 1.38 in Mordukhovich [21], part (i) of this theorem implies that $T$ is Mordukhovich differentiable on $L^2(\Omega)$ such that, for any $h \in L^2(\Omega)$, the Mordukhovich derivative of $T$ at $h$ satisfies that, for any $g \in L^2(\Omega)$,

$$\widehat{D}^*T(h)(g)(t) = T^*g(t)$$

$$= \sum_{j=1}^\infty \left( \sum_{i=1}^\infty c_{ij} \langle g, \varphi_i \rangle \right) \varphi_j(t), \text{ for any } t \in \Omega.$$

Proof of (iii). Since $T$ is a Hilbert-Schmidt integral operator, then, part (iii) can immediately follow from Theorem (3.1). However, in addition to the proof of Theorem 3.1, for the Hilbert-Schmidt integral operator $T$ on $L^2(\Omega)$, we directly prove that

$$\hat{\alpha}(T, f, u) = 0, \text{ for any } f, u \in L^2(\Omega) \text{ with } u = Tf.$$

By definition and by (5.6), for any $m > 1$, we have

$$\hat{\alpha}(T, f, u) = \sup_{\eta > 0} \inf \{ \|v\|_{L^2(\Omega)} : v = \widehat{D}^*T(h)(g), h \in \mathbb{B}(f, \eta), T(h) \in \mathbb{B}(u, \eta), \|g\|_{L^2(\Omega)} = 1 \}$$

$$\leq \sup_{\eta > 0} \inf \{ \|v\|_{L^2(\Omega)} : v = \widehat{D}^*T(h)(\varphi_m), h \in \mathbb{B}(f, \eta), T(h) \in \mathbb{B}(u, \eta), \|g\|_{L^2(\Omega)} = 1, g = \varphi_m \}$$

$$= \sup_{\eta > 0} \inf \left\{ \left\| \sum_{j=1}^\infty \left( \sum_{i=1}^\infty c_{ij} \langle \varphi_m, \varphi_i \rangle \right) \varphi_j \right\|_{L^2(\Omega)} : T^*\varphi_m, h \in \mathbb{B}(f, \eta), T(h) \in \mathbb{B}(u, \eta), g = \varphi_m \right\}$$

$$= \sup_{\eta > 0} \inf \left\{ \left\| \sum_{j=1}^\infty c_{mj} \varphi_j \right\|_{L^2(\Omega)} : T^*\varphi_m, h \in \mathbb{B}(f, \eta), T(h) \in \mathbb{B}(u, \eta), g = \varphi_m \right\}$$

$$= \sup_{\eta>0} \inf \left\{ \sqrt{\sum_{j=1}^{\infty} c_{mj}^2} : \varphi_m, m = 1, 2, \ldots \right\}. \tag{5.7}$$

Since the Hilbert-Schmidt norm of $T$ satisfies

$$\|T^*\|_{HS}^2 = \|k\|_{L^2(\Omega \times \Omega)}^2 = \sum_{i,j=1}^{\infty} c_{ij}^2 < \infty.$$

By $\sum_{i,j=1}^{\infty} c_{ij}^2 = \sum_{m=1}^{\infty} \sum_{j=1}^{\infty} c_{mj}^2 < \infty$, we have

$$\sum_{j=1}^{\infty} c_{mj}^2 \to 0, \text{ as } m \to \infty. \tag{5.8}$$

Applying (5.7) into (5.8), we obtain

$$\hat{\alpha}(T, f, u) = 0, \text{ for any } f, u \in L^2(\Omega) \text{ with } u = Tf. \qquad \square$$

We consider an important class of Hilbert-Schmidt kernel $k$. In particular, suppose that there are $\varphi, \psi \in L^2(\Omega)$ such that the Hilbert-Schmidt kernel $k$ of the considered Hilbert-Schmidt integral operator $T$ studied in this section has the following separable representation

$$k(s, t) = \varphi(s)\psi(t), \text{ for any } (s, t) \in \Omega \times \Omega. \tag{5.9}$$

**Corollary 5.2**. *Suppose $T$ is a Hilbert-Schmidt integral operator with its Hilbert-Schmidt kernel $k$ defined by* (5.9). *Then, $T$ has the following properties.*

(ii)      *$T$ is Fréchet differentiable on $L^2(\Omega)$ such that, for any $h \in L^2(\Omega)$, we have*

$$\nabla T(h)(f) = Tf = \langle f, \psi \rangle \varphi, \text{ for any } \in L^2(\Omega).$$

(ii)      *$T$ is Mordukhovich differentiable on $L^2(\Omega)$ such that, for any $h \in L^2(\Omega)$, we have*

$$\widehat{D}^* T(h)(g) = T^* g = \langle g, \varphi \rangle \psi, \text{ for any } g \in L^2(\Omega).$$

*Proof*. By (5.5) and (5.6), for any $f, g \in L^2(\Omega)$, we have

$$T(f)(s) = \varphi(s) \int_\Omega \psi(t) f(t) dt = \langle f, \psi \rangle \varphi(s), \text{ for any } s \in \Omega.$$

And
$$T^*(g)(t) = \psi(t) \int_\Omega \varphi(s) g(s) ds = \langle g, \varphi \rangle \psi(t), \text{ for any } t \in \Omega.$$

Then this corollary follows from Theorem 5.1 immediately. $\qquad \square$

### 6. Quasi-Hilbert-Schmidt Operators

In this section, we extend the concepts of Hilbert-Schmidt operators on Hilbert spaces to more general cases, in which the linearity is not required. In this section, we also follow the notations used in Section 3. Let $(H, \|\cdot\|_H)$ be a separable Hilbert space with inner product $\langle \cdot, \cdot \rangle$ and equipped with an orthonormal (Schauder) basis $\{e_i : i \in \mathbb{N}\}$ and with origin $\theta$.

**Definition 6.1**. Let $T: H \to H$ be a continuous (the linearity is unrequired) operator. If $T$ satisfies the following condition

$$\sum_{i=1}^{\infty}\|Te_i\|_H^2 < \infty, \tag{6.1}$$

then, $T$ is called a quasi-Hilbert-Schmidt operator on $H$. By the sum in (6.1), the Hilbert-Schmidt norm of $T$ is also defined by

$$\|T\|_{HS}^2 = \sum_{i=1}^{\infty}\|Te_i\|_H^2. \tag{6.2}$$

In this section, we consider a quasi-Hilbert-Schmidt operator on $H$ with finite dimensional range. We will calculate its Fréchet and Mordukhovich derivatives. We define a mapping $T: H \to H$, for any $x \in H$, by

$$Tx = \begin{cases} \frac{\langle x,e_1\rangle^2-\langle x,e_2\rangle^2}{\sqrt{\langle x,e_1\rangle^2+\langle x,e_2\rangle^2}}e_1 + \frac{2\langle x,e_1\rangle\langle x,e_2\rangle}{\sqrt{\langle x,e_1\rangle^2+\langle x,e_2\rangle^2}}e_2, & \text{if } \langle x,e_1\rangle^2 + \langle x,e_2\rangle^2 > 0, \\ \theta, & \text{if } \langle x,e_1\rangle^2 + \langle x,e_2\rangle^2 = 0. \end{cases} \tag{6.3}$$

Associating the mapping $T$ defined by (6.1), for any $z = \sum_{i=1}^{\infty}\langle z,e_i\rangle e_i \in H$, if $\langle z,e_1\rangle^2 + \langle z,e_2\rangle^2 > 0$, then, we write $D(z)$ for the following double sequence of real numbers

$$D(z) = \begin{pmatrix} \frac{(\langle z,e_1\rangle^2+3\langle z,e_2\rangle^2)\langle z,e_1\rangle}{(\langle z,e_1\rangle^2+\langle z,e_2\rangle^2)\sqrt{\langle z,e_1\rangle^2+\langle z,e_2\rangle^2}} & \frac{2\langle z,e_2\rangle^2\langle z,e_2\rangle}{(\langle z,e_1\rangle^2+\langle z,e_2\rangle^2)\sqrt{\langle z,e_1\rangle^2+\langle z,e_2\rangle^2}} & 0 \\ \frac{-(3\langle z,e_1\rangle^2+\langle z,e_2\rangle^2)\langle z,e_2\rangle}{(\langle z,e_1\rangle^2+\langle z,e_2\rangle^2)\sqrt{\langle z,e_1\rangle^2+\langle z,e_2\rangle^2}} & \frac{2\langle z,e_1\rangle^2\langle z,e_1\rangle}{(\langle z,e_1\rangle^2+\langle z,e_2\rangle^2)\sqrt{\langle z,e_1\rangle^2+\langle z,e_2\rangle^2}} & 0 \\ 0 & & 0 \end{pmatrix}. \tag{6.4}$$

Let $D(z)^T$ denote the transpose of $D(z)$. Under the condition that $\langle z,e_1\rangle^2 + \langle z,e_2\rangle^2 > 0$, we have

$$D(z)^T = \begin{pmatrix} \frac{(\langle z,e_1\rangle^2+3\langle z,e_2\rangle^2)\langle z,e_1\rangle}{(\langle z,e_1\rangle^2+\langle z,e_2\rangle^2)\sqrt{\langle z,e_1\rangle^2+\langle z,e_2\rangle^2}} & \frac{-(3\langle z,e_1\rangle^2+\langle z,e_2\rangle^2)\langle z,e_2\rangle}{(\langle z,e_1\rangle^2+\langle z,e_2\rangle^2)\sqrt{\langle z,e_1\rangle^2+\langle z,e_2\rangle^2}} & 0 \\ \frac{2\langle z,e_2\rangle^2\langle z,e_2\rangle}{(\langle z,e_1\rangle^2+\langle z,e_2\rangle^2)\sqrt{\langle z,e_1\rangle^2+\langle z,e_2\rangle^2}} & \frac{2\langle z,e_1\rangle^2\langle z,e_1\rangle}{(\langle z,e_1\rangle^2+\langle z,e_2\rangle^2)\sqrt{\langle z,e_1\rangle^2+\langle z,e_2\rangle^2}} & 0 \\ 0 & & 0 \end{pmatrix}. \tag{6.5}$$

Now, with the notations of $D(z)$ and $D(z)^T$ in (6.4) and (6.5), respectively, we state and prove our main theorem in this section.

**Theorem 6.1**. *Let $(H, \|\cdot\|_H)$ be a separable Hilbert space with dimension great than or equal to* 3. *Let $T: H \to H$ be the operator on $H$ defined in* (6.3). *Then, $T$ has the following properties*:

(i)  $T: H \to H$ *is a nonexpansive mapping*; *that is*

$$\|Tx\|_H \leq \|x\|_H, \text{ for any } x \in H$$

*(ii)*  $T: H \to H$ *is a* (*nonlinear*) *continuous quasi-Hilbert-Schmidt operator on $H$. The Hilbert-Schmidt norm of $T$ satisfies that*

$$\|T\|_{HS}^2 = 2.$$

(iii)  *Let $z = \sum_{i=1}^{\infty}\langle z,e_i\rangle e_i \in H$. Suppose $\langle z,e_1\rangle^2 + \langle z,e_2\rangle^2 > 0$. Then, $T$ is Fréchet differentiable at $z$ such that the Fréchet derivative of $T$ at $z$ satisfies that*

$$\nabla T(z)(x) = \text{sc}(x)D(z)\begin{pmatrix} e_1 \\ e_2 \\ \vdots \end{pmatrix}, \text{ for } x \in H \text{ with } \text{sc}(x) = (\langle x,e_1\rangle,\langle x,e_2\rangle,\ldots) \in l_2.$$

(iv) Let $z = \sum_{i=1}^{\infty} \langle z, e_i \rangle e_i \in H$. Suppose $\langle z, e_1 \rangle^2 + \langle z, e_2 \rangle^2 > 0$. Then, $T$ is Mordukhovich differentiable at $z$ such that,

$$\widehat{D}^*T(z) = \nabla T(z)^*, \text{ this is, } \widehat{D}^*T(z)(y) = \mathrm{sc}(y)D(z)^T \begin{pmatrix} e_1 \\ e_2 \\ \vdots \end{pmatrix}, \text{ for any } y \in H. \tag{6.6}$$

More precisely speaking, for $y = \sum_{i=1}^{\infty} \langle y, e_i \rangle e_i \in H$, (6.6) is represented as

$$\widehat{D}^*T(z)(y) = \left( \langle y, e_1 \rangle \frac{(\langle z,e_1\rangle^2 + 3\langle z,e_2\rangle^2)\langle z,e_1\rangle}{(\langle z,e_1\rangle^2 + \langle z,e_2\rangle^2)\sqrt{\langle z,e_1\rangle^2 + \langle z,e_2\rangle^2}} + \langle y, e_2 \rangle \frac{2\langle z,e_2\rangle^2 \langle z,e_2\rangle}{(\langle z,e_1\rangle^2 + \langle z,e_2\rangle^2)\sqrt{\langle z,e_1\rangle^2 + \langle z,e_2\rangle^2}} \right) e_1$$

$$+ \left( \langle y, e_1 \rangle \frac{-(3\langle z,e_1\rangle^2 + \langle z,e_2\rangle^2)\langle z,e_2\rangle}{(\langle z,e_1\rangle^2 + \langle z,e_2\rangle^2)\sqrt{\langle z,e_1\rangle^2 + \langle z,e_2\rangle^2}} + \langle y, e_2 \rangle \frac{2\langle z,e_1\rangle^2 \langle z,e_1\rangle}{(\langle z,e_1\rangle^2 + \langle z,e_2\rangle^2)\sqrt{\langle z,e_1\rangle^2 + \langle z,e_2\rangle^2}} \right) e_2.$$

(v) Let $z = \sum_{i=1}^{\infty} \langle z, e_i \rangle e_i \in H$. Suppose $\langle z, e_1 \rangle^2 + \langle z, e_2 \rangle^2 > 0$. Then, for any $y \in H$, we have

$$\left\| \widehat{D}^*T(z)(y) \right\|_H^2 = \langle y, e_1 \rangle^2 + \langle y, e_2 \rangle^2 + 12 \left( \langle y, e_1 \rangle \frac{\langle z,e_1\rangle\langle z,e_1\rangle}{\langle z,e_1\rangle^2 + \langle z,e_2\rangle^2} - \langle y, e_1 \rangle \frac{\langle z,e_1\rangle^2 - \langle z,e_2\rangle^2}{2(\langle z,e_1\rangle^2 + \langle z,e_2\rangle^2)} \right)^2.$$

(vi) $T$ maps every bounded subset of $H$ to a relatively compact subset of $H$.

(vii) The covering constant for $T$ satisfies

$$\hat{\alpha}(T, \bar{x}, \bar{y}) = 0, \text{ for any } \bar{x}, \bar{y} \in H \text{ with } \bar{y} = T\bar{x}.$$

(viii) For any $z \in H$, we have

$$\theta \in \widehat{D}^*T(z)(y), \text{ for } y \in H \text{ with } \langle y, e_1 \rangle^2 + \langle y, e_2 \rangle^2 = 0.$$

(ix) Let $z \in H$. If $\langle z, e_1 \rangle^2 + \langle z, e_2 \rangle^2 = 0$, then $T$ is not Mordukhovich differentiable at $z$. More precisely speaking, we have

$$\widehat{D}^*T(z)(y) = \emptyset, \text{ for } y \in H \text{ with } \langle y, e_1 \rangle^2 + \langle y, e_2 \rangle^2 > 0.$$

(x) Let $z \in H$. If $\langle z, e_1 \rangle^2 + \langle z, e_2 \rangle^2 = 0$, then $T$ is not Fréchet differentiable at $z$.

*Proof.* Proof of (i). It is clear that if $\langle x, e_1 \rangle^2 + \langle x, e_2 \rangle^2 = 0$, then $\|Tx\|_H^2 = 0$. Next, suppose $\langle x, e_1 \rangle^2 + \langle x, e_2 \rangle^2 > 0$. We have

$$\|Tx\|_H^2 = \left\| \frac{\langle x,e_1\rangle^2 - \langle x,e_2\rangle^2}{\sqrt{\langle x,e_1\rangle^2 + \langle x,e_2\rangle^2}} e_1 + \frac{2\langle x,e_1\rangle\langle x,e_2\rangle}{\sqrt{\langle x,e_1\rangle^2 + \langle x,e_2\rangle^2}} e_2 \right\|_H^2$$

$$= \left( \frac{\langle x,e_1\rangle^2 - \langle x,e_2\rangle^2}{\sqrt{\langle x,e_1\rangle^2 + \langle x,e_2\rangle^2}} \right)^2 + \left( \frac{2\langle x,e_1\rangle\langle x,e_2\rangle}{\sqrt{\langle x,e_1\rangle^2 + \langle x,e_2\rangle^2}} \right)^2$$

$$= \langle x, e_1 \rangle^2 + \langle x, e_2 \rangle^2$$

$$\leq \|x\|_H^2.$$

This proves (i). Then, we prove (ii). In particular, in definition (6.3), let $x = e_k$, for $k = 1, 2, \ldots$, we have

$$Te_k = \frac{\langle e_k, e_1\rangle^2 - \langle e_k, e_2\rangle^2}{\sqrt{\langle e_k, e_1\rangle^2 + \langle e_k, e_2\rangle^2}} e_1 + \frac{2\langle e_k, e_1\rangle \langle e_k, e_2\rangle}{\sqrt{\langle e_k, e_1\rangle^2 + \langle e_k, e_2\rangle^2}} e_2.$$

This implies that $Te_1 = e_1$, $Te_2 = -e_1$ and $Te_k = \theta$, for $k > 2$. Then, by (6.2), we calculate the Hilbert-Schmidt norm of $T$.

$$\|T\|_{HS}^2 = \sum_{k=1}^{\infty} \|Te_k\|_H^2 = \|Te_1\|_H^2 + \|Te_2\|_H^2 = \|e_1\|_H^2 + \|-e_1\|_H^2 = 2.$$

This proves (ii). Then, we prove (iii). In [1, 20, 23], a mapping $f: \mathbb{R}^2 \to \mathbb{R}^2$ is defined by

$$f(x) = f(x_1, x_2) = \left(\frac{x_1^2 - x_2^2}{\sqrt{x_1^2 + x_2^2}}, \frac{2x_1 x_2}{\sqrt{x_1^2 + x_2^2}}\right), \text{ for } x = (x_1, x_2) \in \mathbb{R}^2 \setminus \{\theta\} \text{ with } f(\theta) = \theta. \quad (6.7)$$

In [20], the Fréchet and Mordukhovich differentiability of this mapping $f$ defined by (6.7) were studied, by which its covering constant was calculated. We recall some results about this mapping $f$ from [20].

**Theorem 5.3 in [20].** *Let $z = (z_1, z_2) \in \mathbb{R}^2 \setminus \{\theta\}$. Then $f$ is Fréchet differentiable and Mordukhovich differentiable at $z$ such that*

$$\nabla f(z) = \begin{pmatrix} \frac{(z_1^2 + 3z_2^2)z_1}{(z_1^2 + z_2^2)\sqrt{z_1^2 + z_2^2}} & \frac{2z_2^2 z_2}{(z_1^2 + z_2^2)\sqrt{z_1^2 + z_2^2}} \\ \frac{-(3z_1^2 + z_2^2)z_2}{(z_1^2 + z_2^2)\sqrt{z_1^2 + z_2^2}} & \frac{2z_1^2 z_1}{(z_1^2 + z_2^2)\sqrt{z_1^2 + z_2^2}} \end{pmatrix} \quad \text{and} \quad \widehat{D}^* f(z) = \begin{pmatrix} \frac{(z_1^2 + 3z_2^2)z_1}{(z_1^2 + z_2^2)\sqrt{z_1^2 + z_2^2}} & \frac{-(3z_1^2 + z_2^2)z_2}{(z_1^2 + z_2^2)\sqrt{z_1^2 + z_2^2}} \\ \frac{2z_2^2 z_2}{(z_1^2 + z_2^2)\sqrt{z_1^2 + z_2^2}} & \frac{2z_1^2 z_1}{(z_1^2 + z_2^2)\sqrt{z_1^2 + z_2^2}} \end{pmatrix}.$$

*However, more precisely speaking, for $x = (x_1, x_2), y = (y_1, y_2) \in \mathbb{R}^2$, if $x = \widehat{D}^* f(z)(y)$, then*

$$x_1 = y_1 \frac{(z_1^2 + 3z_2^2)z_1}{(z_1^2 + z_2^2)\sqrt{z_1^2 + z_2^2}} + y_2 \frac{2z_2^2 z_2}{(z_1^2 + z_2^2)\sqrt{z_1^2 + z_2^2}},$$

$$x_2 = y_1 \frac{(-3z_1^2 - z_2^2)z_2}{(z_1^2 + z_2^2)\sqrt{z_1^2 + z_2^2}} + y_2 \frac{2z_1^2 z_1}{(z_1^2 + z_2^2)\sqrt{z_1^2 + z_2^2}}.$$

Let $\|\cdot\|_2$ denote the Hilbert norm in $\mathbb{R}^2$. By the definition of Fréchet derivatives, we have

$$\lim_{x \to z} \frac{f(x) - f(z) - \nabla f(z)(x - z)}{\|x - z\|_2} = (0, 0).$$

This is

$$\lim_{x \to z} \frac{\left(\frac{x_1^2 - x_2^2}{\sqrt{x_1^2 + x_2^2}}, \frac{2x_1 x_2}{\sqrt{x_1^2 + x_2^2}}\right) - \left(\frac{z_1^2 - z_2^2}{\sqrt{z_1^2 + z_2^2}}, \frac{2z_1 z_2}{\sqrt{z_1^2 + z_2^2}}\right) - (x_1 - z_1, x_2 - z_2) \begin{pmatrix} \frac{(z_1^2 + 3z_2^2)z_1}{(z_1^2 + z_2^2)\sqrt{z_1^2 + z_2^2}} & \frac{2z_2^2 z_2}{(z_1^2 + z_2^2)\sqrt{z_1^2 + z_2^2}} \\ \frac{-(3z_1^2 + z_2^2)z_2}{(z_1^2 + z_2^2)\sqrt{z_1^2 + z_2^2}} & \frac{2z_1^2 z_1}{(z_1^2 + z_2^2)\sqrt{z_1^2 + z_2^2}} \end{pmatrix}}{\|x - z\|_2} = (0, 0). \quad (6.8)$$

More precisely speaking, by considering each coordinate, (6.8) implies that

$$\lim_{x \to z} \frac{\frac{x_1^2 - x_2^2}{\sqrt{x_1^2 + x_2^2}} - \frac{z_1^2 - z_2^2}{\sqrt{z_1^2 + z_2^2}} - \left((x_1 - z_1)\frac{(z_1^2 + 3z_2^2)z_1}{(z_1^2 + z_2^2)\sqrt{z_1^2 + z_2^2}} + (x_2 - z_2)\frac{-(3z_1^2 + z_2^2)z_2}{(z_1^2 + z_2^2)\sqrt{z_1^2 + z_2^2}}\right)}{\|x - z\|_2} = 0, \tag{6.9}$$

and

$$\lim_{x \to z} \frac{\frac{2x_1 x_2}{\sqrt{x_1^2 + x_2^2}} - \frac{2z_1 z_2}{\sqrt{z_1^2 + z_2^2}} - \left((x_1 - z_1)\frac{2z_2^2 z_2}{(z_1^2 + z_2^2)\sqrt{z_1^2 + z_2^2}} + (x_2 - z_2)\frac{2z_1^2 z_1}{(z_1^2 + z_2^2)\sqrt{z_1^2 + z_2^2}}\right)}{\|x - z\|_2} = 0. \tag{6.10}$$

In (6.9) and (6.10), when $x = (x_1, x_2)$ and $z = (z_1, z_2)$ are simultaneously and respectively replaced by $\hat{x} = (\langle x, e_1 \rangle, \langle x, e_2 \rangle)$ and $\hat{z} = (\langle z, e_1 \rangle, \langle z, e_2 \rangle)$, we obtain that

$$\lim_{\hat{x} \to \hat{z}} \frac{\frac{\langle x, e_1 \rangle^2 - \langle x, e_2 \rangle^2}{\sqrt{\langle x, e_1 \rangle^2 + \langle x, e_2 \rangle^2}} - \frac{\langle z, e_1 \rangle^2 - \langle z, e_2 \rangle^2}{\sqrt{\langle z, e_1 \rangle^2 + \langle z, e_2 \rangle^2}} - \left((\langle x, e_1 \rangle - \langle z, e_1 \rangle)\frac{(\langle z, e_1 \rangle^2 + 3\langle z, e_2 \rangle^2)\langle z, e_1 \rangle}{(\langle z, e_1 \rangle^2 + \langle z, e_2 \rangle^2)\sqrt{\langle z, e_1 \rangle^2 + \langle z, e_2 \rangle^2}} + (\langle x, e_2 \rangle - \langle z, e_2 \rangle)\frac{-(3\langle z, e_1 \rangle^2 + \langle z, e_2 \rangle^2)\langle z, e_2 \rangle}{(\langle z, e_1 \rangle^2 + \langle z, e_2 \rangle^2)\sqrt{\langle z, e_1 \rangle^2 + \langle z, e_2 \rangle^2}}\right)}{\|\hat{x} - \hat{z}\|_2} = 0, \tag{6.11}$$

$$\lim_{\hat{x} \to \hat{z}} \frac{\frac{2\langle x, e_1 \rangle \langle x, e_2 \rangle}{\sqrt{\langle x, e_1 \rangle^2 + \langle x, e_2 \rangle^2}} - \frac{2\langle z, e_1 \rangle \langle z, e_2 \rangle}{\sqrt{\langle z, e_1 \rangle^2 + \langle z, e_2 \rangle^2}} - \left((\langle x, e_1 \rangle - \langle z, e_1 \rangle)\frac{2\langle z, e_2 \rangle^2 \langle z, e_2 \rangle}{(\langle z, e_1 \rangle^2 + \langle z, e_2 \rangle^2)\sqrt{\langle z, e_1 \rangle^2 + \langle z, e_2 \rangle^2}} + (\langle x, e_2 \rangle - \langle z, e_2 \rangle)\frac{2\langle z, e_1 \rangle^2 \langle z, e_1 \rangle}{(\langle z, e_1 \rangle^2 + \langle z, e_2 \rangle^2)\sqrt{\langle z, e_1 \rangle^2 + \langle z, e_2 \rangle^2}}\right)}{\|\hat{x} - \hat{z}\|_2} = 0. \tag{6.12}$$

By combining (6.11) and (6.12), we get

$$\lim_{\hat{x} \to \hat{z}} \frac{\left(\frac{\langle x, e_1 \rangle^2 - \langle x, e_2 \rangle^2}{\sqrt{\langle x, e_1 \rangle^2 + \langle x, e_2 \rangle^2}} - \frac{\langle z, e_1 \rangle^2 - \langle z, e_2 \rangle^2}{\sqrt{\langle z, e_1 \rangle^2 + \langle z, e_2 \rangle^2}}\right) e_1 - \left((\langle x, e_1 \rangle - \langle z, e_1 \rangle)\frac{(\langle z, e_1 \rangle^2 + 3\langle z, e_2 \rangle^2)\langle z, e_1 \rangle}{(\langle z, e_1 \rangle^2 + \langle z, e_2 \rangle^2)\sqrt{\langle z, e_1 \rangle^2 + \langle z, e_2 \rangle^2}} + (\langle x, e_2 \rangle - \langle z, e_2 \rangle)\frac{-(3\langle z, e_1 \rangle^2 + \langle z, e_2 \rangle^2)\langle z, e_2 \rangle}{(\langle z, e_1 \rangle^2 + \langle z, e_2 \rangle^2)\sqrt{\langle z, e_1 \rangle^2 + \langle z, e_2 \rangle^2}}\right) e_1}{\|\hat{x} - \hat{z}\|_2}$$
$$= \theta. \tag{6.13}$$

$$\lim_{\hat{x} \to \hat{z}} \frac{\left(\frac{2\langle x, e_1 \rangle \langle x, e_2 \rangle}{\sqrt{\langle x, e_1 \rangle^2 + \langle x, e_2 \rangle^2}} - \frac{2\langle z, e_1 \rangle \langle z, e_2 \rangle}{\sqrt{\langle z, e_1 \rangle^2 + \langle z, e_2 \rangle^2}}\right) e_2 - \left((\langle x, e_1 \rangle - \langle z, e_1 \rangle)\frac{2\langle z, e_2 \rangle^2 \langle z, e_2 \rangle}{(\langle z, e_1 \rangle^2 + \langle z, e_2 \rangle^2)\sqrt{\langle z, e_1 \rangle^2 + \langle z, e_2 \rangle^2}} + (\langle x, e_2 \rangle - \langle z, e_2 \rangle)\frac{2\langle z, e_1 \rangle^2 \langle z, e_1 \rangle}{(\langle z, e_1 \rangle^2 + \langle z, e_2 \rangle^2)\sqrt{\langle z, e_1 \rangle^2 + \langle z, e_2 \rangle^2}}\right) e_2}{\|\hat{x} - \hat{z}\|_2}$$
$$= \theta. \tag{6.14}$$

For $z = \sum_{i=1}^{\infty} \langle z, e_i \rangle e_i \in H$ satisfying $\langle z, e_1 \rangle^2 + \langle z, e_2 \rangle^2 > 0$, by the definition of $T$ in (6.3), the limits (6.13) and (6.14) can be equivalently rewritten as

$$\lim_{\hat{x} \to \hat{z}} \frac{Tx - Tz - (\hat{x} - \hat{z})\begin{pmatrix} \frac{(\langle z, e_1 \rangle^2 + 3\langle z, e_2 \rangle^2)\langle z, e_1 \rangle}{(\langle z, e_1 \rangle^2 + \langle z, e_2 \rangle^2)\sqrt{\langle z, e_1 \rangle^2 + \langle z, e_2 \rangle^2}} & \frac{2\langle z, e_2 \rangle^2 \langle z, e_2 \rangle}{(\langle z, e_1 \rangle^2 + \langle z, e_2 \rangle^2)\sqrt{\langle z, e_1 \rangle^2 + \langle z, e_2 \rangle^2}} \\ \frac{-(3\langle z, e_1 \rangle^2 + \langle z, e_2 \rangle^2)\langle z, e_2 \rangle}{(\langle z, e_1 \rangle^2 + \langle z, e_2 \rangle^2)\sqrt{\langle z, e_1 \rangle^2 + \langle z, e_2 \rangle^2}} & \frac{2\langle z, e_1 \rangle^2 \langle z, e_1 \rangle}{(\langle z, e_1 \rangle^2 + \langle z, e_2 \rangle^2)\sqrt{\langle z, e_1 \rangle^2 + \langle z, e_2 \rangle^2}} \end{pmatrix}\begin{pmatrix} e_1 \\ e_2 \end{pmatrix}}{\|\hat{x} - \hat{z}\|_2} = \theta. \tag{6.15}$$

For $x, z \in H$, it is clearly to see that $\|\hat{x} - \hat{z}\|_2 \leq \|x - z\|_H$. Then, (6.15) implies that

$$\lim_{x \to z} \frac{Tx - Tz - (\hat{x}-\hat{z}) \begin{pmatrix} \frac{(\langle z,e_1 \rangle^2 + 3\langle z,e_2 \rangle^2)\langle z,e_1 \rangle}{(\langle z,e_1 \rangle^2 + \langle z,e_2 \rangle^2)\sqrt{\langle z,e_1 \rangle^2 + \langle z,e_2 \rangle^2}} & \frac{2\langle z,e_2 \rangle^2 \langle z,e_2 \rangle}{(\langle z,e_1 \rangle^2 + \langle z,e_2 \rangle^2)\sqrt{\langle z,e_1 \rangle^2 + \langle z,e_2 \rangle^2}} \\ \frac{-(3\langle z,e_1 \rangle^2 + \langle z,e_2 \rangle^2)\langle z,e_2 \rangle}{(\langle z,e_1 \rangle^2 + \langle z,e_2 \rangle^2)\sqrt{\langle z,e_1 \rangle^2 + \langle z,e_2 \rangle^2}} & \frac{2\langle z,e_1 \rangle^2 \langle z,e_1 \rangle}{(\langle z,e_1 \rangle^2 + \langle z,e_2 \rangle^2)\sqrt{\langle z,e_1 \rangle^2 + \langle z,e_2 \rangle^2}} \end{pmatrix} \begin{pmatrix} e_1 \\ e_2 \end{pmatrix}}{\|x-z\|_H} = \theta. \quad (6.16)$$

By the definition of the double sequence of real numbers $D(z)$ in (6.4), for $x = \sum_{i=1}^{\infty} \langle x, e_i \rangle e_i$, and $z = \sum_{i=1}^{\infty} \langle z, e_i \rangle e_i \in H$ satisfying $\langle z, e_1 \rangle^2 + \langle z, e_2 \rangle^2 > 0$, (6.16) can be rewritten as

$$\lim_{x \to z} \frac{Tx - Tz - ((\langle x,e_1 \rangle, \langle x,e_2 \rangle, \dots) - (\langle z,e_1 \rangle, \langle z,e_2 \rangle, \dots))D(z)\begin{pmatrix} e_1 \\ e_2 \\ \vdots \end{pmatrix}}{\|x-z\|_H} = \theta. \quad (6.17)$$

Recall that, for $x = \sum_{i=1}^{\infty} \langle x, e_i \rangle e_i \in H$, $\text{sc}(x) = (\langle x, e_1 \rangle, \langle x, e_2 \rangle, \dots) \in l_2$. Then, (6.17) is equivalent to

$$\lim_{x \to z} \frac{Tx - Tz - (\text{sc}(x) - \text{sc}(z))D(z)\begin{pmatrix} e_1 \\ e_2 \\ \vdots \end{pmatrix}}{\|x-z\|_H} = \theta. \quad (6.18)$$

(6.18) implies that, for $z = \sum_{i=1}^{\infty} \langle z, e_i \rangle e_i \in H$ satisfying $\langle z, e_1 \rangle^2 + \langle z, e_2 \rangle^2 > 0$, $T$ is Fréchet differentiable at $z$ such that the Fréchet derivative of $T$, $\nabla T(z)$, is a linear and continuous mapping from $H$ to $H$ such that

$$\nabla T(z)(x) = \text{sc}(x)D(z)\begin{pmatrix} e_1 \\ e_2 \\ \vdots \end{pmatrix}, \text{ for any } x \in H. \quad (6.19)$$

This proves (iii). Then, we prove (iv). Let $z = \sum_{i=1}^{\infty} \langle z, e_i \rangle e_i \in H$ satisfying $\langle z, e_1 \rangle^2 + \langle z, e_2 \rangle^2 > 0$. By (6.19) and by Theorem 1.38 in Mordukhovich [21], $T$ is Mordukhovich differentiable at $z$ such that

$$\widehat{D}^*T(z) = \nabla T(z)^*.$$

Here, $\nabla T(z)^*$ is the adjoint operator of $\nabla T(z)$. Then, $\nabla T(z)^*$ is defined by $D(z)^T$, which is the transpose of $D(z)$ given in (6.5). That is,

$$\widehat{D}^*T(z)(y) = \text{sc}(y)D(z)^T \begin{pmatrix} e_1 \\ e_2 \\ \vdots \end{pmatrix}, \text{ for any } y \in H. \quad (6.20)$$

More precisely speaking, for $y = \sum_{i=1}^{\infty} \langle y, e_i \rangle e_i \in H$, (6.20) is represented as

$\widehat{D}^*T(z)(y)$

$$= (\langle y, e_1 \rangle, \langle y, e_2 \rangle, \dots) \begin{pmatrix} \frac{(\langle z,e_1 \rangle^2 + 3\langle z,e_2 \rangle^2)\langle z,e_1 \rangle}{(\langle z,e_1 \rangle^2 + \langle z,e_2 \rangle^2)\sqrt{\langle z,e_1 \rangle^2 + \langle z,e_2 \rangle^2}} & \frac{-(3\langle z,e_1 \rangle^2 + \langle z,e_2 \rangle^2)\langle z,e_2 \rangle}{(\langle z,e_1 \rangle^2 + \langle z,e_2 \rangle^2)\sqrt{\langle z,e_1 \rangle^2 + \langle z,e_2 \rangle^2}} & 0 \\ \frac{2\langle z,e_2 \rangle^2 \langle z,e_2 \rangle}{(\langle z,e_1 \rangle^2 + \langle z,e_2 \rangle^2)\sqrt{\langle z,e_1 \rangle^2 + \langle z,e_2 \rangle^2}} & \frac{2\langle z,e_1 \rangle^2 \langle z,e_1 \rangle}{(\langle z,e_1 \rangle^2 + \langle z,e_2 \rangle^2)\sqrt{\langle z,e_1 \rangle^2 + \langle z,e_2 \rangle^2}} & 0 \\ 0 & 0 \end{pmatrix} \begin{pmatrix} e_1 \\ e_2 \\ \vdots \end{pmatrix}$$

$$= \left( \langle y, e_1 \rangle \frac{(\langle z,e_1 \rangle^2 + 3\langle z,e_2 \rangle^2)\langle z,e_1 \rangle}{(\langle z,e_1 \rangle^2 + \langle z,e_2 \rangle^2)\sqrt{\langle z,e_1 \rangle^2 + \langle z,e_2 \rangle^2}} + \langle y, e_2 \rangle \frac{2\langle z,e_2 \rangle^2 \langle z,e_2 \rangle}{(\langle z,e_1 \rangle^2 + \langle z,e_2 \rangle^2)\sqrt{\langle z,e_1 \rangle^2 + \langle z,e_2 \rangle^2}} \right) e_1$$

$$+ \left( \langle y, e_1 \rangle \frac{-(3\langle z,e_1\rangle^2+\langle z,e_2\rangle^2)\langle z,e_2\rangle}{(\langle z,e_1\rangle^2+\langle z,e_2\rangle^2)\sqrt{\langle z,e_1\rangle^2+\langle z,e_2\rangle^2}} + \langle y, e_2 \rangle \frac{2\langle z,e_1\rangle^2\langle z,e_1\rangle}{(\langle z,e_1\rangle^2+\langle z,e_2\rangle^2)\sqrt{\langle z,e_1\rangle^2+\langle z,e_2\rangle^2}} \right) e_2. \quad (6.21)$$

Proof of (v). By (6.21), we have

$$\left\| \widehat{D}^* T(z)(y) \right\|_H^2 = \left( \langle y, e_1 \rangle \frac{(\langle z,e_1\rangle^2+3\langle z,e_2\rangle^2)\langle z,e_1\rangle}{(\langle z,e_1\rangle^2+\langle z,e_2\rangle^2)\sqrt{\langle z,e_1\rangle^2+\langle z,e_2\rangle^2}} + \langle y, e_2 \rangle \frac{2\langle z,e_2\rangle^2\langle z,e_2\rangle}{(\langle z,e_1\rangle^2+\langle z,e_2\rangle^2)\sqrt{\langle z,e_1\rangle^2+\langle z,e_2\rangle^2}} \right)^2$$

$$+ \left( \langle y, e_1 \rangle \frac{-(3\langle z,e_1\rangle^2+\langle z,e_2\rangle^2)\langle z,e_2\rangle}{(\langle z,e_1\rangle^2+\langle z,e_2\rangle^2)\sqrt{\langle z,e_1\rangle^2+\langle z,e_2\rangle^2}} + \langle y, e_2 \rangle \frac{2\langle z,e_1\rangle^2\langle z,e_1\rangle}{(\langle z,e_1\rangle^2+\langle z,e_2\rangle^2)\sqrt{\langle z,e_1\rangle^2+\langle z,e_2\rangle^2}} \right)^2$$

$$= \langle y, e_1 \rangle^2 + \langle y, e_2 \rangle^2 + 12 \left( \langle y, e_1 \rangle \frac{\langle z,e_1\rangle\langle z,e_1\rangle}{\langle z,e_1\rangle^2+\langle z,e_2\rangle^2} - \langle y, e_1 \rangle \frac{\langle z,e_1\rangle^2-\langle z,e_2\rangle^2}{2(\langle z,e_1\rangle^2+\langle z,e_2\rangle^2)} \right)^2.$$

The proof of the last equation is similar to the proof of Proposition 5.4 in [20] and it is omitted here.

Proof of (vi). Since the range of $T$ has dimension 2, then, part (vi) is proved immediately.

Proof of (vii). Part (vii) is an immediate consequence of part (vi) and Theorem 3.2. However, we give a direct proof for part (vii). Let $\bar{x}, \bar{y} \in H$ with $\bar{y} = T\bar{x}$. By Parts (iv) and (v), we calculate the covering constant for $T$ at point $(\bar{x}, \bar{y})$. Since $(H, \|\cdot\|_H)$ is a separable Hilbert space with dimension great than or equal to 3, we have

$$\hat{\alpha}(T, \bar{x}, \bar{y}) = \sup_{\eta>0} \inf \{ \|x\|_H : x = \widehat{D}^* T(z)(y), z \in \mathbb{B}(\bar{x}, \eta), T(z) \in \mathbb{B}(\bar{y}, \eta), \|y\|_H = 1 \}$$

$$\leq \sup_{\eta>0} \inf \{ \|x\|_H : x = \widehat{D}^* T(z)(y), z \in \mathbb{B}(\bar{x}, \eta), \langle z, e_1 \rangle^2 + \langle z, e_2 \rangle^2 > 0, T(z) \in \mathbb{B}(\bar{y}, \eta), \|y\|_H = 1 \}$$

$$= \sup_{\eta>0} \inf \left\{ \sqrt{\langle y, e_1 \rangle^2 + \langle y, e_2 \rangle^2 + 12 \left( \langle y, e_1 \rangle \frac{\langle z,e_1\rangle\langle z,e_1\rangle}{\langle z,e_1\rangle^2+\langle z,e_2\rangle^2} - \langle y, e_1 \rangle \frac{\langle z,e_1\rangle^2-\langle z,e_2\rangle^2}{2(\langle z,e_1\rangle^2+\langle z,e_2\rangle^2)} \right)^2 } : \right.$$

$$z \in \mathbb{B}(\bar{x}, \eta), \langle z, e_1 \rangle^2 + \langle z, e_2 \rangle^2 > 0, T(z) \in \mathbb{B}(\bar{y}, \eta), \|y\|_H = 1 \}$$

$$\leq \sup_{\eta>0} \inf \left\{ \sqrt{\langle y, e_1 \rangle^2 + \langle y, e_2 \rangle^2 + 12 \left( \langle y, e_1 \rangle \frac{\langle z,e_1\rangle\langle z,e_1\rangle}{\langle z,e_1\rangle^2+\langle z,e_2\rangle^2} - \langle y, e_1 \rangle \frac{\langle z,e_1\rangle^2-\langle z,e_2\rangle^2}{2(\langle z,e_1\rangle^2+\langle z,e_2\rangle^2)} \right)^2 } : \right.$$

$$z \in \mathbb{B}(\bar{x}, \eta), \langle z, e_1 \rangle^2 + \langle z, e_2 \rangle^2 > 0, T(z) \in \mathbb{B}(\bar{y}, \eta), \|y\|_H = 1, \langle y, e_1 \rangle = \langle y, e_2 \rangle = 0 \}$$

$$= 0.$$

Notice that, the last equality follows from the assumption in this theorem that $(H, \|\cdot\|_H)$ is a separable Hilbert space with dimension great than or equal to 3.

In the following proof of this theorem, we adopt some notations. For $x = \sum_{i=1}^{\infty} \langle x, e_i \rangle e_i \in H$, we write $x_i = \langle x, e_i \rangle$, for $i = 1, 2, \ldots$.

Proof of (viii). Let $z$ be an arbitrarily given point in $H$. Let $y \in H$ with $\langle y, e_1 \rangle^2 + \langle y, e_2 \rangle^2 = 0$. Notice that by the definition (6.1) and (6.3), that $T$ maps $H$ into a two dimensional subspace of $H$, which is generated by $\{e_1, e_2\}$. It implies that $\langle T(u) - T(z), e_i \rangle = 0$, for $i = 3, 4, \ldots$. We have

$$\limsup_{u \to z} \frac{\langle \theta, u-z \rangle - \langle y, T(u)-T(z) \rangle}{\|u-z\| + \|T(u)-T(z)\|}$$

$$= \limsup_{u \to z} \frac{0 - \langle (0,0,y_3,y_4,\ldots),\ T(u)-T(z)\rangle}{\|u-z\|+\|T(u)-T(z)\|}$$

$$= 0.$$

This proves (viii). Then, we prove (ix). Let $z = \sum_{i=1}^{\infty}\langle z, e_i\rangle e_i \in H$. Suppose that $\langle z, e_1\rangle^2 + \langle z, e_2\rangle^2 = 0$. Let $y \in H$ with $\langle y, e_1\rangle^2 + \langle y, e_2\rangle^2 > 0$. By the condition that $z_1 = \langle z, e_1\rangle = z_2 = \langle z, e_2\rangle = 0$ and by the definition of $T$, we have that $Tz = \theta$ and, for any $x \in H$, $\langle Tu, e_i\rangle = 0$, for $i = 3, 4, \ldots$. For $x \in H$, we calculate

$$\limsup_{u \to z} \frac{\langle x,\ u-z\rangle - \langle y,\ T(u)-T(z)\rangle}{\|u-z\|+\|T(u)-T(z)\|}$$

$$= \limsup_{u \to z} \frac{\sum_{i=1}^{\infty} x_i(u_i-z_i) - \left(y_1\left(\frac{u_1^2-u_2^2}{\sqrt{u_1^2+u_2^2}}-0\right) + y_2\left(\frac{2u_1 u_2}{\sqrt{u_1^2+u_2^2}}-0\right)\right)}{\|u-z\|+\|T(u)-T(z)\|}$$

$$= \limsup_{u \to z} \frac{x_1(u_1-0)+x_2(u_2-0) - y_1\left(\frac{u_1^2-u_2^2}{\sqrt{u_1^2+u_2^2}}-0\right) - y_2\left(\frac{2u_1 u_2}{\sqrt{u_1^2+u_2^2}}-0\right) + \sum_{i=3}^{\infty} x_i(u_i-z_i)}{\|u-z\|+\|T(u)-(0,0)\|}$$

$$= \limsup_{u \to z} \frac{x_1 u_1 + x_2 u_2 - y_1 \frac{u_1^2-u_2^2}{\sqrt{u_1^2+u_2^2}} - y_2 \frac{2u_1 u_2}{\sqrt{u_1^2+u_2^2}} + \sum_{i=3}^{\infty} x_i(u_i-z_i)}{\sqrt{u_1^2+u_2^2+\sum_{i=3}^{\infty}(u_i-z_i)^2} + \sqrt{\left(\frac{u_1^2-u_2^2}{\sqrt{u_1^2+u_2^2}}\right)^2 + \left(\frac{2u_1 u_2}{\sqrt{u_1^2+u_2^2}}\right)^2}}. \tag{6.22}$$

Case 1. Assume that $x = \theta$. In this case, the limit (6.22) becomes

$$\limsup_{u \to z} \frac{x_1 u_1 + x_2 u_2 - y_1 \frac{u_1^2-u_2^2}{\sqrt{u_1^2+u_2^2}} - y_2 \frac{2u_1 u_2}{\sqrt{u_1^2+u_2^2}} + \sum_{i=3}^{\infty} x_i(u_i-z_i)}{\sqrt{u_1^2+u_2^2+\sum_{i=3}^{\infty}(u_i-z_i)^2} + \sqrt{\left(\frac{u_1^2-u_2^2}{\sqrt{u_1^2+u_2^2}}\right)^2 + \left(\frac{2u_1 u_2}{\sqrt{u_1^2+u_2^2}}\right)^2}}$$

$$= \limsup_{u \to z} \frac{-y_1 \frac{u_1^2-u_2^2}{\sqrt{u_1^2+u_2^2}} - y_2 \frac{2u_1 u_2}{\sqrt{u_1^2+u_2^2}}}{\sqrt{u_1^2+u_2^2+\sum_{i=3}^{\infty}(u_i-z_i)^2} + \sqrt{\left(\frac{u_1^2-u_2^2}{\sqrt{u_1^2+u_2^2}}\right)^2 + \left(\frac{2u_1 u_2}{\sqrt{u_1^2+u_2^2}}\right)^2}}. \tag{6.23}$$

In this case, by the assumption that $\langle y, e_1\rangle^2 + \langle y, e_2\rangle^2 > 0$, we have $y_1 \neq 0$, or $y_2 \neq 0$, or both.

Subcase 1.1. Suppose that $y_1 < 0$. By the condition $z_1 = z_2 = 0$, in this case, we take a special direction in the limit (6.23) for $u \to z$ by $u_2 = 0$ and $u_1 \downarrow 0$ and $u_i = z_i$, for all $i \geq 3$. The limit (6.23) become that

$$\limsup_{u \to z} \frac{-y_1 \frac{u_1^2 - u_2^2}{\sqrt{u_1^2 + u_2^2}} - y_2 \frac{2u_1 u_2}{\sqrt{u_1^2 + u_2^2}}}{\sqrt{u_1^2 + u_2^2 + \sum_{i=3}^{\infty}(u_i - z_i)^2} + \sqrt{\left(\frac{u_1^2 - u_2^2}{\sqrt{u_1^2 + u_2^2}}\right)^2 + \left(\frac{2u_1 u_2}{\sqrt{u_1^2 + u_2^2}}\right)^2}}$$

$$\geq \limsup_{u_1 \downarrow 0} \frac{-y_1 \frac{u_1^2}{\sqrt{u_1^2}}}{\sqrt{u_1^2} + \sqrt{\left(\frac{u_1^2}{\sqrt{u_1^2}}\right)^2}}$$

$$= -\frac{y_1}{2} > 0.$$

This proves that

$$\theta \notin \widehat{D}^* T(z)(y), \text{ for } y \in H \setminus \{\theta\} \text{ with } y_1 < 0. \tag{6.24}$$

It can be similarly proved that

$$\theta \notin \widehat{D}^* T(z)(y), \text{ for } y \in H \setminus \{\theta\} \text{ with } y_1 > 0. \tag{6.25}$$

Subcase 1.2. Suppose that $y_2 < 0$. By the condition $z_1 = z_2 = 0$, in this case, we take a special direction in the limit (6.23) for $u \to z$ by $u_2 = u_1$ and $u_i = z_i$, for all $i \geq 3$ with $u_1 \downarrow 0$. The limit (6.23) become that

$$\limsup_{u \to z} \frac{-y_1 \frac{u_1^2 - u_2^2}{\sqrt{u_1^2 + u_2^2}} - y_2 \frac{2u_1 u_2}{\sqrt{u_1^2 + u_2^2}}}{\sqrt{u_1^2 + u_2^2 + \sum_{i=3}^{\infty}(u_i - z_i)^2} + \sqrt{\left(\frac{u_1^2 - u_2^2}{\sqrt{u_1^2 + u_2^2}}\right)^2 + \left(\frac{2u_1 u_2}{\sqrt{u_1^2 + u_2^2}}\right)^2}}$$

$$\geq \limsup_{u_1 \downarrow 0} \frac{-y_2 \frac{2u_1^2}{\sqrt{2u_1^2}}}{\sqrt{2u_1^2} + \sqrt{\left(\frac{2u_1^2}{\sqrt{2u_1^2}}\right)^2}}$$

$$= -\frac{y_2}{2} > 0.$$

This proves that

$$\theta \notin \widehat{D}^* T(z)(y), \text{ for } y \in H \setminus \{\theta\} \text{ with } y_2 < 0. \tag{6.26}$$

It can be similarly proved that

$$\theta \notin \widehat{D}^* T(z)(y), \text{ for } y \in H \setminus \{\theta\} \text{ with } y_2 > 0. \tag{6.27}$$

By (2.24) − (2.27), we obtain that

$$\theta \notin \widehat{D}^* T(z)(y), \text{ for } y \in H \setminus \{\theta\} \text{ with } \langle y, e_1 \rangle^2 + \langle y, e_2 \rangle^2 > 0. \tag{6.28}$$

Case 2. Let $x \in H\setminus\{\theta\}$. Case 2 is considered by the following subcases.

Subcase 2.1. Assume that $\sum_{i=3}^{\infty} x_i^2 > 0$, that is there are some $i \geq 3$ such that $x_i \neq 0$. In this case, we take a special direction in the limit (6.22) for $u \to z$ by $u_1 = z_1 = u_2 = z_2 = 0$ and $u_i - z_i = \frac{1}{n}x_i$, for all $i \geq 3$, with $n \to \infty$. By the assumption that $\sum_{i=3}^{\infty} x_i^2 > 0$, the limit (6.22) becomes that

$$\limsup_{u \to z} \frac{x_1 u_1 + x_2 u_2 - y_1 \frac{u_1^2 - u_2^2}{\sqrt{u_1^2 + u_2^2}} - y_2 \frac{2u_1 u_2}{\sqrt{u_1^2 + u_2^2}} + \sum_{i=3}^{\infty} x_i(u_i - z_i)}{\sqrt{u_1^2 + u_2^2 + \sum_{i=3}^{\infty}(u_i - z_i)^2} + \sqrt{\left(\frac{u_1^2 - u_2^2}{\sqrt{u_1^2 + u_2^2}}\right)^2 + \left(\frac{2u_1 u_2}{\sqrt{u_1^2 + u_2^2}}\right)^2}}$$

$$\geq \limsup_{u \to z} \frac{x_1 \cdot 0 + x_2 \cdot 0 - y_1 \cdot 0 - y_2 \cdot 0 + \sum_{i=3}^{\infty} x_i(u_i - z_i)}{\sqrt{0 + 0 + \sum_{i=3}^{\infty}(u_i - z_i)^2} + 0}$$

$$\geq \limsup_{n \to \infty} \frac{\sum_{i=3}^{\infty} \frac{1}{n} x_i^2}{\sqrt{\sum_{i=3}^{\infty} \frac{1}{n^2} x_i^2}}$$

$$= \limsup_{n \to \infty} \sqrt{\sum_{i=3}^{\infty} x_i^2}$$

$$> 0.$$

This implies that

$$x \notin \widehat{D}^*T(z)(y), \text{ for } x \in H\setminus\{\theta\} \text{ with } \sum_{i=3}^{\infty} x_i^2 > 0. \tag{6.29}$$

Subcase 2.2. Assume that $\sum_{i=3}^{\infty} x_i^2 = 0$, that is $x_i = 0$, for all $i \geq 3$, and $x_1^2 + x_2^2 > 0$. In this case, the limit (6.22) becomes

$$\limsup_{u \to z} \frac{x_1 u_1 + x_2 u_2 - y_1 \frac{u_1^2 - u_2^2}{\sqrt{u_1^2 + u_2^2}} - y_2 \frac{2u_1 u_2}{\sqrt{u_1^2 + u_2^2}} + \sum_{i=3}^{\infty} x_i(u_i - z_i)}{\sqrt{u_1^2 + u_2^2 + \sum_{i=3}^{\infty}(u_i - z_i)^2} + \sqrt{\left(\frac{u_1^2 - u_2^2}{\sqrt{u_1^2 + u_2^2}}\right)^2 + \left(\frac{2u_1 u_2}{\sqrt{u_1^2 + u_2^2}}\right)^2}}$$

$$= \limsup_{u \to z} \frac{x_1 u_1 + x_2 u_2 - y_1 \frac{u_1^2 - u_2^2}{\sqrt{u_1^2 + u_2^2}} - y_2 \frac{2u_1 u_2}{\sqrt{u_1^2 + u_2^2}} + \sum_{i=3}^{\infty} 0 \cdot (u_i - z_i)}{\sqrt{u_1^2 + u_2^2 + \sum_{i=3}^{\infty}(u_i - z_i)^2} + \sqrt{\left(\frac{u_1^2 - u_2^2}{\sqrt{u_1^2 + u_2^2}}\right)^2 + \left(\frac{2u_1 u_2}{\sqrt{u_1^2 + u_2^2}}\right)^2}}$$

$$= \limsup_{u \to z} \frac{x_1 u_1 + x_2 u_2 - y_1 \frac{u_1^2 - u_2^2}{\sqrt{u_1^2 + u_2^2}} - y_2 \frac{2u_1 u_2}{\sqrt{u_1^2 + u_2^2}}}{\sqrt{u_1^2 + u_2^2 + \sum_{i=3}^{\infty}(u_i - z_i)^2} + \sqrt{\left(\frac{u_1^2 - u_2^2}{\sqrt{u_1^2 + u_2^2}}\right)^2 + \left(\frac{2u_1 u_2}{\sqrt{u_1^2 + u_2^2}}\right)^2}}. \tag{6.30}$$

In the limit (6.30), we take a special direction for $u \to z$ as $u_i = z_i$, for all $i \geq 3$. Notice that $z_1 = z_2 =$

0. Then, this special direction is equivalent to $(u_1, u_2) \to (z_1, z_2)$ in $\mathbb{R}^2$. With respect to this special direction, limit (6.30) satisfies that

$$\limsup_{u \to z} \frac{x_1 u_1 + x_2 u_2 - y_1 \frac{u_1^2 - u_2^2}{\sqrt{u_1^2 + u_2^2}} - y_2 \frac{2u_1 u_2}{\sqrt{u_1^2 + u_2^2}}}{\sqrt{u_1^2 + u_2^2 + \sum_{i=3}^{\infty}(u_i - z_i)^2} + \sqrt{\left(\frac{u_1^2 - u_2^2}{\sqrt{u_1^2 + u_2^2}}\right)^2 + \left(\frac{2u_1 u_2}{\sqrt{u_1^2 + u_2^2}}\right)^2}}$$

$$\geq \limsup_{u \to z} \frac{x_1 u_1 + x_2 u_2 - y_1 \frac{u_1^2 - u_2^2}{\sqrt{u_1^2 + u_2^2}} - y_2 \frac{2u_1 u_2}{\sqrt{u_1^2 + u_2^2}}}{\sqrt{u_1^2 + u_2^2} + \sqrt{\left(\frac{u_1^2 - u_2^2}{\sqrt{u_1^2 + u_2^2}}\right)^2 + \left(\frac{2u_1 u_2}{\sqrt{u_1^2 + u_2^2}}\right)^2}} \quad \text{(with respect to this special direction)}$$

$$= \limsup_{(u_1, u_2) \to (z_1, z_2) \text{ in } \mathbb{R}^2} \frac{x_1 u_1 + x_2 u_2 - y_1 \frac{u_1^2 - u_2^2}{\sqrt{u_1^2 + u_2^2}} - y_2 \frac{2u_1 u_2}{\sqrt{u_1^2 + u_2^2}}}{\sqrt{u_1^2 + u_2^2} + \sqrt{\left(\frac{u_1^2 - u_2^2}{\sqrt{u_1^2 + u_2^2}}\right)^2 + \left(\frac{2u_1 u_2}{\sqrt{u_1^2 + u_2^2}}\right)^2}}$$

$$= \limsup_{(u_1, u_2) \to (0,0) \text{ in } \mathbb{R}^2} \frac{x_1 u_1 + x_2 u_2 - y_1 \frac{u_1^2 - u_2^2}{\sqrt{u_1^2 + u_2^2}} - y_2 \frac{2u_1 u_2}{\sqrt{u_1^2 + u_2^2}}}{\sqrt{u_1^2 + u_2^2} + \sqrt{\left(\frac{u_1^2 - u_2^2}{\sqrt{u_1^2 + u_2^2}}\right)^2 + \left(\frac{2u_1 u_2}{\sqrt{u_1^2 + u_2^2}}\right)^2}}$$

$$= \limsup_{(u_1, u_2) \to (0,0) \text{ in } \mathbb{R}^2} \frac{x_1 u_1 + x_2 u_2 - y_1 \frac{u_1^2 - u_2^2}{\sqrt{u_1^2 + u_2^2}} - y_2 \frac{2u_1 u_2}{\sqrt{u_1^2 + u_2^2}}}{2\sqrt{u_1^2 + u_2^2}}. \tag{6.31}$$

For the assume in Subcase 2.2 that $\sum_{i=3}^{\infty} x_i^2 = 0$, that is $x_i = 0$, for all $i \geq 3$ and $x_1^2 + x_2^2 = \langle x, e_1 \rangle^2 + \langle x, e_2 \rangle^2 > 0$. Meanwhile, the limit (6.31) is considered under the condition that $y \in H \setminus \{\theta\}$ with $y_1^2 + y_2^2 = \langle y, e_1 \rangle^2 + \langle y, e_2 \rangle^2 > 0$. Then, the estimation of limit (6.31) is divided into the following subcase, which are studied in [20]. For details, see the proof of Proposition 5.2 in [20].

Subcase 2.2.1. $x_1 > y_1$. In this case, it was proved in [20] that

$$\limsup_{u \to \theta} \frac{x_1 u_1 + x_2 u_2 - y_1 \frac{u_1^2 - u_2^2}{\sqrt{u_1^2 + u_2^2}} - y_2 \frac{2u_1 u_2}{\sqrt{u_1^2 + u_2^2}}}{2\|u\|} \geq \frac{x_1}{2} - \frac{y_1}{2} > 0. \tag{6.32}$$

Subcase 2.2.2. $-x_1 > y_1$. In this case, it was proved in [20] that

$$\limsup_{u \to \theta} \frac{x_1 u_1 + x_2 u_2 - y_1 \frac{u_1^2 - u_2^2}{\sqrt{u_1^2 + u_2^2}} - y_2 \frac{2u_1 u_2}{\sqrt{u_1^2 + u_2^2}}}{2\|u\|} \geq -\frac{x_1}{2} - \frac{y_1}{2} > 0. \tag{6.33}$$

Subcase 2.2.3. $x_1 + x_2 > \sqrt{2} y_2$. In this case, it was proved in [20] that

$$\limsup_{u \to \theta} \frac{x_1 u_1 + x_2 u_2 - y_1 \frac{u_1^2 - u_2^2}{\sqrt{u_1^2 + u_2^2}} - y_2 \frac{2 u_1 u_2}{\sqrt{u_1^2 + u_2^2}}}{2 \|u\|} = \frac{x_1 + x_2 - \sqrt{2} y_2}{2\sqrt{2}} > 0. \tag{6.34}$$

Subcase 2.2.4. $-x_1 - x_2 > \sqrt{2} y_2$. In this case, it was proved in [20] that

$$\limsup_{u \to \theta} \frac{x_1 u_1 + x_2 u_2 - y_1 \frac{u_1^2 - u_2^2}{\sqrt{u_1^2 + u_2^2}} - y_2 \frac{2 u_1 u_2}{\sqrt{u_1^2 + u_2^2}}}{2 \|u\|} \geq \frac{-x_1 - x_2 - \sqrt{2} y_2}{2\sqrt{2}} > 0. \tag{6.35}$$

Subcase 2.2.5. $-x_1 + x_2 > -\sqrt{2} y_2$. In this case, it was proved in [20] that

$$\limsup_{u \to \theta} \frac{x_1 u_1 + x_2 u_2 - y_1 \frac{u_1^2 - u_2^2}{\sqrt{u_1^2 + u_2^2}} - y_2 \frac{2 u_1 u_2}{\sqrt{u_1^2 + u_2^2}}}{2 \|u\|} \geq \frac{-x_1}{2\sqrt{2}} + \frac{x_2}{2\sqrt{2}} + \frac{\sqrt{2} y_2}{2\sqrt{2}} > 0. \tag{6.36}$$

Subcase 2.2.6. $x_1 - x_2 > -\sqrt{2} y_2$. In this case, it was proved in [20] that

$$\limsup_{u \to \theta} \frac{x_1 u_1 + x_2 u_2 - y_1 \frac{u_1^2 - u_2^2}{\sqrt{u_1^2 + u_2^2}} - y_2 \frac{2 u_1 u_2}{\sqrt{u_1^2 + u_2^2}}}{2 \|u\|} \geq \frac{x_1}{2\sqrt{2}} - \frac{x_2}{2\sqrt{2}} + \frac{\sqrt{2} y_2}{2\sqrt{2}} > 0. \tag{6.37}$$

The inequalities (6.32) to (6.37) imply that, we obtain that if $x$ satisfies any one of the conditions listed in subcases 2.2.1 to 2.2.6, then

$$\limsup_{(u_1, u_2) \to (0,0) \text{ in } \mathbb{R}^2} \frac{x_1 u_1 + x_2 u_2 - y_1 \frac{u_1^2 - u_2^2}{\sqrt{u_1^2 + u_2^2}} - y_2 \frac{2 u_1 u_2}{\sqrt{u_1^2 + u_2^2}}}{2 \sqrt{u_1^2 + u_2^2}} > 0. \tag{6.38}$$

(6.38) implies that

$x \notin \widehat{D}^* T(z)(y)$, if $x$ satisfies any one of the conditions listed in subcases 2.2.1 to 2.2.6. (6.39)

By (6.38) we have that for $x \in H \setminus \{\theta\}$ to satisfy the following inequality

$$\limsup_{(u_1, u_2) \to (0,0) \text{ in } \mathbb{R}^2} \frac{x_1 u_1 + x_2 u_2 - y_1 \frac{u_1^2 - u_2^2}{\sqrt{u_1^2 + u_2^2}} - y_2 \frac{2 u_1 u_2}{\sqrt{u_1^2 + u_2^2}}}{2 \sqrt{u_1^2 + u_2^2}} \leq 0,$$

the following inequalities are necessary conditions (opposite conditions of subcases 2.2.1 to 2.2.6)

(i) $x_1 \leq y_1$, (ii) $-x_1 \leq y_1$, (iii) $x_1 + x_2 \leq \sqrt{2} y_2$,

(iv) $-x_1 - x_2 \leq \sqrt{2} y_2$, (v) $-x_1 + x_2 \leq -\sqrt{2} y_2$, (vi) $x_1 - x_2 \leq -\sqrt{2} y_2$.

By (6.39), the above inequalities (i−vi) are necessary conditions for $x \in \widehat{D}^* f(\theta)(y)$. However, in [20], it has been proved that the unique solution of the above inequalities (i−vi) for $x_1, x_2$ is

$$x_1 = x_2 = 0.$$

This contradicts to the assumption that $x_1^2 + x_2^2 = \langle x, e_1 \rangle^2 + \langle x, e_2 \rangle^2 > 0$, by which, we obtain that

$$x \notin \widehat{D}^*T(z)(y), \text{ if } x \text{ satisfies inequalities (i–vi) and } x_1^2 + x_2^2 > 0. \tag{6.40}$$

By combining (6.39) and (6.40), we have

$$x \notin \widehat{D}^*T(z)(y), \text{ for } x \in H\setminus\{\theta\} \text{ with } \langle x, e_1 \rangle^2 + \langle x, e_2 \rangle^2 > 0. \tag{6.41}$$

This proves the subcase 2.2 that

$$x \notin \widehat{D}^*T(z)(y), \text{ for } x \in H\setminus\{\theta\} \text{ with } \sum_{i=3}^{\infty} x_i^2 = 0. \tag{6.42}$$

By (4.29) and (4.42), we prove case 2 that, for $y \in H\setminus\{\theta\}$ with $\langle y, e_1 \rangle^2 + \langle y, e_2 \rangle^2 > 0$, we have

$$x \notin \widehat{D}^*T(z)(y), \text{ for } x \in H\setminus\{\theta\}. \tag{6.43}$$

By (6.28) in case 1 and (6.43) in case 2, we obtain that

$$\widehat{D}^*T(z)(y) = \emptyset, \text{ for } y \in H\setminus\{\theta\} \text{ with } \langle y, e_1 \rangle^2 + \langle y, e_2 \rangle^2 > 0.$$

This proves (ix). Part (x) follows from Theorem 1.38 in [21] and part (ix) immediately. □

## 7. Conclusion and Remarks

It is known that the generalized differentiation theory in Banach spaces is based on Fréchet derivatives of single-valued mappings and Mordukhovich derivatives (coderivatives) of both single-valued and set-valued mappings. Recall that in the theory of (classical) differentiation in calculus, there are many "derivative formulas" for real valued functions. These formulas firmly formed the foundation of applications of differentiation. In contrast to the differentiation theory in calculus, in order to have some real application in generalized differentiation theory, we must find some "Fréchet derivative formulas" and some "Mordukhovich derivative formulas". More precisely speaking, we want to find the exact Fréchet derivatives of some single-valued mappings and find the exact Mordukhovich derivatives of some single-valued mappings.

For this purpose, in this paper, we calculate the Fréchet derivatives and Mordukhovich derivatives of Hilbert-Schmidt operators on separable Hilbert spaces, which includes Hilbert-Schmidt integral operators, as an important class. We introduce the concept of quasi-Hilbert-Schmidt operators on separable Hilbert spaces and find their Fréchet derivatives and Mordukhovich derivatives.

It is important to notice that in this paper, the Mordukhovich derivatives of Hilbert-Schmidt operators and quasi-Hilbert-Schmidt operators on separable Hilbert spaces are calculated by using the connection between Fréchet and Mordukhovich derivatives, which is Theorem 1.38 in [21]. It raises some questions for our consideration.

(Q1) For a given single-valued mapping on Hilbert spaces (on Banach spaces, in general), how do we directly and precisely find its Mordukhovich derivative without using Theorem 1.38 in [21]?

The following question should be more difficult to answer.

(Q2) For a given set-valued mapping on Euclidean spaces (on Hilbert spaces, or on Banach spaces, in general), how do we precisely find its Mordukhovich derivative?

# References

[1] Arutyunov A. V., Mordukhovich B. S. and Zhukovskiy S. E., Coincidence Points of Parameterized Generalized Equations with Applications to Optimal Value Functions, Journal of Optimization Theory and Applications 196, 177–198 (2023).

[2] Arutyunov, A.V., Avakov, E.R., Zhukovskiy, S.E.: Stability theorems for estimating the distance to a set of coincidence points. SIAM J. Optim. 25, 807–828 (2015).

[3] Arutyunov A. V. and S. E. Zhukovskiy, Hadamard's theorem for mappings with relaxed smoothness conditions, Sbornik: Mathematics,210, 165–183 (2019).

[4] Arutyunov, A.V., Zhukovskiy, S.E.: Existence and properties of inverse mappings. Proc. Steklov Inst. Math. 271, 12–22 (2010).

[5] Arutyunov, A.V.: Smooth abnormal problems in extremum theory and analysis. Russ. Math. Surv. 67, 403–457 (2012).

[6] Bao, T.Q., Gupta, P., Mordukhovich, B.S.: Necessary conditions in multiobjective optimization with equilibrium constraints. J. Optim. Theory Appl. 135, 179–203 (2007).

[7] Bonnans, J. F., Shapiro, A.: Perturbation Analysis of Optimization Problems. Springer, New York (2000).

[8] Bump, Daniel. *Automorphic Forms and Representations*. Cambridge University Press. ISBN 0-521-65818-7 (1998).

[9] Coleman, Rodney, ed. Calculus on Normed Vector Spaces, Universitext, Springer, ISBN 978-1-4614-3894-6. (2012).

[10] Conway, John B. A course in functional analysis. New York: Springer-Verlag. ISBN 978-0-387-97245-9. OCLC 21195908 (1990).

[11] Dieudonné, Jean *Foundations of modern analysis, Boston, MA:* Academic Press, MR 0349288 (1969).

[12] Dontchev, A.L., Rockefeller, R.T.: Implicit Functions and Solution Mappings. A View from Variational Analysis, Springer, New York (2014).

[13] Lang, Serge *Differential and Riemannian Manifolds,* Springer, ISBN 0-387-94338-2 (1995).

[14] Li, J. L., Mordukhovich derivatives of the metric projection operator in Hilbert spaces, *Journal of Optimization Theory and Applications* (2024).

[15] Li, J. L., Mordukhovich derivatives of the metric projection operator in uniformly convex and uniformly smooth Banach spaces, *Set-Valued and Variational Analysis* (2024).

[16] Li, J. L., Covering Constants for Metric Projection Operator with Applications to Stochastic Fixed-Point Problems, *Journal of Global Optimization* (2025).

[17] Li, J. L., Some Applications of Arutyunov Mordukhovich and Zhukovskiy Theorem to Stochastic Integral Equations, to appear in *Journal of Nonlinear and Convex Analysis*, arXiv:2511.03623.

[18] Li, J. L., Petrusel, A., Wen, C F., and Yao, J. C., Covering Constants for Linear and Continuous Mappings with Applications to Stochastic Systems of Linear Equations, to appear in *Journal of Nonlinear and Variational Analysis.*

[19] Li, J. L., Tammer C., and Yao, J. C., Some Applications of Arutyunov Mordukhovich Zhukovskiy Theorem to Set-Valued Stochastic Vector Variational Inequalities and Optimizations, in preparation.

[20] Li, J. L., Calculating Covering Constants for Mappings in Euclidean Spaces Using Mordukhovich Derivatives (Coderivatives) with Applications, submitted.

[21] Mordukhovich, B.S.: Complete characterizations of openness, metric regularity, and Lipschitzian properties of multifunctions. Trans. Amer. Math. Soc. 340, 1–35 (1993).

[22] Mordukhovich, B.S.: Variational Analysis and Applications. Springer, Switzerland (2018).

[23] Mordukhovich, B.S.: Variational Analysis and Generalized Differentiation, I: Basic Theory, II: Applications. Springer, Berlin (2006).

[24] Mordukhovich, B.S., Nam, N.M.: Convex Analysis and Beyond, I: Basic Theory. Springer, Cham, Switzerland (2022).

[25] Moslehian, M. S., Hilbert–Schmidt Operator (From MathWorld).

[26] Munkres, James R. *Analysis on manifolds,* Addison-Wesley, ISBN 978-0-201-51035-5, MR 1079066. (1991).

[27] Reich, S., A remark on a problem of Asplund, Atti Accad. Lincei 67, 204—205 (1979).

[28] Renardy, Michael; Rogers, Robert C. *An Introduction to Partial Differential Equations.* New York Berlin Heidelberg: Springer Science & Business Media. ISBN 0-387-00444-0 (2004).

[29] Rockafeller, R.T., Wets, R.J.-B.: Variational Analysis. Springer, Berlin (1998).

[30]  Schaefer, Helmut H., *Topological Vector Spaces*. GTM. Vol. 3. New York, NY: Springer New York Imprint Springer. ISBN 978-1-4612-7155-0. OCLC 840278135 (1999).
[31]  Simon, B., *An Overview of Rigorous Scattering Theory*. S2CID 16913591(1978).
[32]  Voitsekhovskii, M. I., Hilbert-Schmidt operator, Encyclopedia of Mathematics, EMS Press (2001).